\documentclass[10pt,a4paper]{report}
\usepackage{wnewthesis,graphicx}
\usepackage{ epsfig}
\usepackage[centertags]{amsmath}
\usepackage{amsfonts}
\usepackage{amsthm}
\usepackage{amssymb,  amscd, subfigure, fancyhdr}

\usepackage[total={6in,9in}, left=1in, right=1in, top=1in, bottom=1in, includefoot]{geometry}

\newtheorem{The}{Theorem}[chapter]
\newtheorem{Lem}{Lemma}[chapter]
\newtheorem{Cor}{Corollary}[chapter]

\newtheorem{Prop}{Proposition}[chapter]
\newtheorem*{Folklore}{Folklore}
\newtheorem*{GCT}{Geodesic Characterization Theorem}
\newtheorem*{StThe}{Stallings' Theorem}

\newtheorem*{Conc4}{Condition $C^{\prime \prime} (4)$}
\newtheorem*{Cont4}{Condition $T (4)$ }
\newtheorem*{Geocomthe}{Geodesic Completion Theorem}

\theoremstyle{remark}
\newtheorem{Remark}{Remark}[chapter]

\theoremstyle{definition}
\newtheorem{Def}{Definition}[chapter]
\newtheorem{Alg}{Algorithm}[chapter]
\newtheorem*{Duwo}{Algorithm \ref{dugoword}$^{\prime}$}

 \leftchapter                       

\singlespace                       

\renewcommand{\thesisdepartmentname}{Mathematics Institute}    

\renewcommand{\thesissubmission}{Submitted to the University of Warwick\\
                        for the degree of}

\renewcommand{\thesisauthor}{Iain Moffatt}    

\renewcommand{\thesismonth}{November}     

\renewcommand{\thesisyear}{2004}      

\renewcommand{\thesistitle}{Integration and Conjugacy in Knot Theory}     

\begin{document}

\thesiscopyrightpage                 


 \thesistitlecolourpage           

\tableofcontents                     
 \listoffigures                     

\begin{thesisacknowledgments}        
  First and foremost I would like to thank my supervisor Stavros Garoufalidis for his  support, guidance and suggestions.
Most of this thesis was written while a visitor at Georgia Institute of Technology between November 2003 and May 2004 and I would like to thank the math department for its hospitality, EPSRC for providing the bulk of the financial support and Colin Rourke for being my Warwick contact.
Finally I would like thank Daan Krammer, whose suggestions greatly improved chapter 2, Werner Nickel for some GAP code and Daniel Groves, Derek Holt, Chuck Miller and Alvaro Pelayo for sharing some of their mathematical knowledge with me.                    


\end{thesisacknowledgments}

\begin{thesisdeclaration}

Except where otherwise stated, this thesis is my own work.
I confirm that this thesis has not been submitted for a degree at any other university. Some of the results presented here have been submitted for  publication.

\end{thesisdeclaration}          

\begin{thesisabstract}
   This thesis consists of three self-contained chapters. The first two concern quantum invariants of links and three manifolds and the third contains results on the word problem for link groups.

In chapter~1 we relate the tree part of the \AA rhus integral to the $\mu$-invariants of string-links in homology balls thus generalizing  results of Habegger and Masbaum.

There is a folklore result in physics saying that the Feynman integration of an exponential is itself an exponential. In chapter~2 we state and prove an exact formulation of this statement in the language which is used in the theory of finite type invariants.

The final chapter is concerned with properties of link groups. In particular we study the relationship between known  solutions from small cancellation theory and normal surface theory for the word and conjugacy problems of the groups of (prime) alternating links. We show that two of the algorithms in the literature for solving the word problem, each using one of the two approaches, are the same.  Then, by considering small cancellation methods, we give a normal surface solution to the conjugacy problem  of these link groups and characterize the conjugacy classes.
Finally as an application of the small cancellation properties of link groups we give a new proof that alternating links are non-trivial.

\end{thesisabstract}

\newpage{\pagestyle{empty}\cleardoublepage} 

\chapter{The $\mu$-invariants and the \AA rhus Integral} \label{chapter:mu}


\def\STU#1#2{      
\begin{array}{ccccc}
\epsfig{file=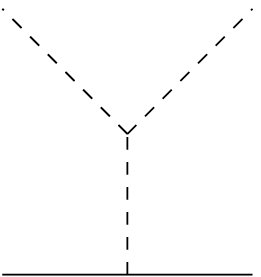, height=#1 } & \raisebox{#2}{=} &\epsfig{file=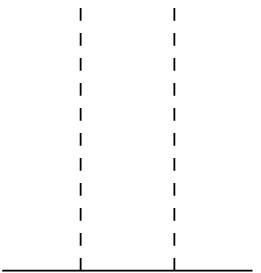, height=#1 } & \raisebox{#2}{-} & \epsfig{file=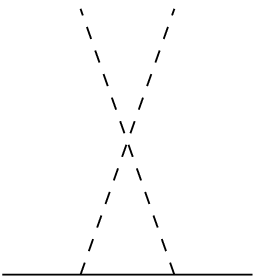, height=#1 }
\end{array}}

\def\IHX#1#2{      
\begin{array}{ccccc}
\epsfig{file=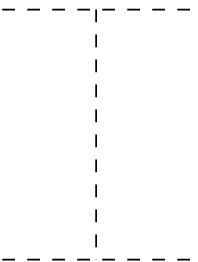, height=#1 } & \raisebox{#2}{=} &\epsfig{file=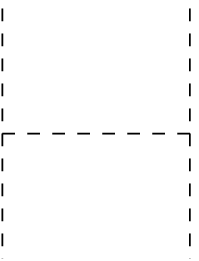, height=#1 } & \raisebox{#2}{-} & \epsfig{file=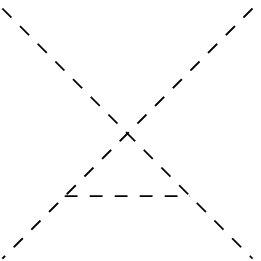, height=#1 }
\end{array}}

\def\AS#1#2{      
\begin{array}{ccc}
\epsfig{file=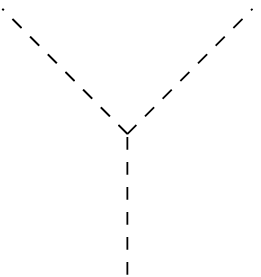, height=#1 } &  \raisebox{#2}{=} & \raisebox{#2}{-} \epsfig{file=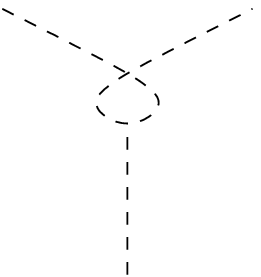, height=#1 }
\end{array}}


\def\bracx#1#2{\langle #1 ,  #2 \rangle_{X} }
\def\FG{\sideset{}{_{X_M}^{FG}}\int}
\def\fgi{\int_{X_M}^{FG}}
\def\bracxm#1#2{\langle #1 ,  #2 \rangle_{X_M} }
\def\bracxmt#1#2{\langle #1 ,  #2 \rangle_{X_M}^{t} }
\def\brat#1#2{\langle #1 ,  #2 \rangle^{t}  }
\def\bracy#1#2{\langle #1 ,  #2 \rangle_{Y} }
\def\bracs#1#2{\langle #1 ,  #2 \rangle_{S} }


\def\paamlt{Z_{0}^{M:t} (\sml )}
\def\aamlt{Z^{M:t} (\sml )}
\def\aart#1{Z^{M;t} (#1)}
\def\paart#1{Z_{0}^{M;t} (#1)}
\def\aar#1{Z^M (#1)}
\def\aaml{Z^M (\sml )}
\def\zm#1{Z^M (#1)}
\def\paaml{Z_{0}^{M} (\sml)}
\def\paar#1{Z_{0}^{M} (#1)}
\def\ea{Z^{M;h,t}(\lambda_i \otimes 1_M)}
\def\eb{\pi^{h,t} (Z^M (1_0 \otimes \sml )^{-1} D_i(Z^M ( \sml ))^{b_i})}
\def\Zmht{Z^{M;h,t}}
\def\zmht#1{Z^{M;h,t} (#1)}
\def\bra#1#2#3{\langle #1 ,  #2 \rangle_{#3} }
\def\pzm#1{Z_{0}^{M} (#1)}
\def\z{\check{Z}}
\def\zmhtn#1{Z_{\leq n}^{M;h,t} (#1)}
\def\Zhtn#1{Z_{\leq n}^{h,t} (#1)}
\def\zht#1{\check{Z}^{h,t} (#1)}
\def\pzm{Z^{M}_{0}}
\def\Zm{Z^{M}}
\def\zm#1{Z^{M} (#1)}
\def\zmt#1{Z^{M;t} (#1)}
\def\ZM{Z^{M}}
\def\Zht#1{Z^{h,t} (#1)}
\def\ZMH{Z^{M;h}}
\def\PZMT{Z_{0}^{M;t}}
\def\ZMT{Z^{M;t}}
\def\ZM{Z^{M}}
\def\PZM{Z^{M}_{0}}

\def\sm{\sigma_{M}}
\def\pa#1{\partial ( #1 )}
\def\bsig{B^{\Sigma}}
\def\h#1#2{H_{#1}(#2)}
\def\fun#1{\pi_1 (#1)}
\def\pres#1#2{\langle #1 | #2 \rangle}
\def\map#1#2#3{#1 : #2 \rightarrow #3}
\def\quoq#1{#1 / #1_{q+1} }
\def\fquo#1{\frac{#1}{#1_{q+1}} }
\def\art{\textit{Art}}
\def\pow#1{\mathcal{P} \textit{(#1)} }
\def\ser#1#2{ #1_1 , \ldots , #1_{#2}}
\def\MU#1#2{\mu_{i}^{(n)}}
\def\cp{\vartriangle}
\def\Exp#1{\exp (#1)}
\def\coeff#1#2{\text{Coeff}(#1,#2)}
\def\prodab#1#2{\prod_{i= #1}^{#2} a_{i}^{b_i}}
\def\ab#1{a_{#1}^{b_{#1}}}
\def\lmi{\lambda_{i}^{M}}
\def\li{\lambda_{i}}
\def\Exp#1{\text{exp} (#1)}
\def\pht{\pi^{h,t}}
\def\sml{\sigma_{LM}}
\def\dti{D^2 \times I}
\def\udot{\mathaccent\cdot\cup}
\def\sxm{\sigma_{X_M}}
\def\mubar{\bar{\mu}}
\def\chiso{\raisebox{1mm}{$\chi$}}
\def\quo#1{#1 / #1_{n+1} }


\def\Cht#1{\mathcal{C}^{h,t} (#1) }
\def\Ct#1{\mathcal{C}^{t} (#1) }

\def\Alm{\mathcal{A}(\uparrow_{X_L \cup X_M})}
\def\alglm{\mathcal{A}(\uparrow_{X_L }, X_M)}
\def\algle{\mathcal{A}(\uparrow_{X_L }, \emptyset )}
\def\Al1m{\mathcal{A}(\uparrow_{X_L +1 \cup X_M})}
\def\algl1e{\mathcal{A}(\uparrow_{X_L+1 }, \emptyset)}
\def\alghtl1e{\mathcal{A}^{h,t}(\uparrow_{X_L+1 }, \emptyset)}
\def\algyz{\mathcal{A}(\uparrow_X, Y)}
\def\algspecial{\mathcal{A}^{h,t}_{n-[ \frac{n}{2} ]-1}(\uparrow_{Y_0 \cup Y_L }, \emptyset)}
\def\Ahtlo{\mathcal{A}^{h,t}(\uparrow_{X_L +1 })}
\def\Bhtlo{\mathcal{B}^{h,t}(X_L +1 )}
\def\algtl1e{\mathcal{A}^t(\uparrow_{X_L+1 }, \emptyset)}
\def\Ahtl1m{\mathcal{A}^{h,t}(\uparrow_{X_L +1 \cup X_M})}
\def\Atl{\mathcal{A}^t(\uparrow_{X_L })}
\def\algtle{\mathcal{A}^{t}(\uparrow_{X_L }, \emptyset)}
\def\algez{\mathcal{A}(\emptyset, Y)}
\def\algyzn{\mathcal{A}_n(\uparrow_{X}, Y)}
\def\algchiimage{\mathcal{A}(\uparrow_{X \cup Y^{\prime}}, Y-Y^{\prime})}
\def\algyzh{\mathcal{A}^h(\uparrow_{X}, Y)}
\def\algyzt{\mathcal{A}^t(\uparrow_{X}, Y)}
\def\algrs{\mathcal{A}(\uparrow_X, Y)}
\def\algqq#1#2{\mathcal{A}(\uparrow_{#1 }, #2)}
\def\algxy{\mathcal{A}(\uparrow_X, Y)}
\def\algxyt{\mathcal{A}^t(\uparrow_X, Y)}
\def\algxyh{\mathcal{A}^h(\uparrow_X, Y)}
\def\Ahoxl1{ \mathcal{A}^{h(0)}(\uparrow_{X_L +1})}
\def\Axl1{ \mathcal{A}(\uparrow_{X_L +1})}
\def\algtlm{\mathcal{A}^{t}(\uparrow_{X_L }, X_M)}

In this chapter we relate the tree part of the \AA rhus integral to the $\mu$-invariants of string-links in homology balls.


\section{Introduction}
Milnor's $\mubar$-invariants of Links and their well defined cousins, the $\mu$-invariants of string-links are classical and well-studied invariants.
These invariants have been brought into the realm of finite-type invariants by Bar-Natan in \cite{BN:95}, Lin in \cite{Li:97} and Habegger and Masbaum in \cite{HM:00}.
Here we are particularly interested in Habegger and Masbaum's formula which expresses the $\mu$-invariants
 in terms of the tree part of the Kontsevich integral.

  The literature on the $\mubar$-invariants is mostly concerned with links in $S^3$.  The generalization to $\mubar$-invariants of links in integral homology spheres and  $\mu$-invariants of string-links in homology balls exists mostly as folklore.  We discuss the $\mu$-invariants of string-links in homology balls and generalize Habegger and Masbaum's results by relating the $\mu$-invariants to the \AA rhus integral, which is a generalization of the Kontsevich integral to links in rational homology spheres, defined in \cite{BGRT:AI, BGRT:AII, BGRT:AIII}.
We do this by representing string-links in homology spheres by string-links in $\dti$ with some distinguished surgery components.

The reader may find some familiarity with the basic properties of the Kontsevich integral useful.

\section{Tangles and String-links} \label{sec:defs}

Let $B^M$ be a connected, compact orientable 3-manifold equipped with a fixed identification $\varphi$ of  the boundary with $\pa{\dti}$.
A {\em tangle} of $n$ components $T \subset B^M$ is a smooth compact 1-manifold $X$, of $n$ components, together with a smooth  embedding
$T: (X, \partial (X) ) \rightarrow (B^M , \partial (B^M))$, transverse to the boundary.
As is standard, we abuse notation and  confuse a tangle, its embedding and its isotopy class.

By a \textit{framing} on a component $i$  of a tangle we mean that we equip  $i$ with a non-vanishing vector field such that the restriction to the boundary is the restriction of a fixed unit vector field normal to the $x$-axis of $D^2$ under the identification.

A {\em coloured} tangle is a tangle equipped with a bijection from the components onto a set of cardinality $n$, where $n$ is the number of components of the tangle.

Since $B^M$ can be obtained by surgery on a framed link $L \subset \dti$, we may represent a tangle $T \subset B^M $ by a tangle $T^{\prime} \subset \dti$ some of whose components are distinguished framed copies of $S^1$, on which we do the surgery. We say that $T^{\prime}$ {\em represents} $T$.
We will call  these distinguished components the {\em surgery components} and the other components the {\em linking components}.
If the tangle $T$ is coloured then this partitions  the colouring set into  sets  corresponding to the surgery components and the linking components. In this chapter we denote these  sets $X_M$ and $X_L$ respectively.

\begin{figure}
\centering
\subfigure[]{\epsfig{file=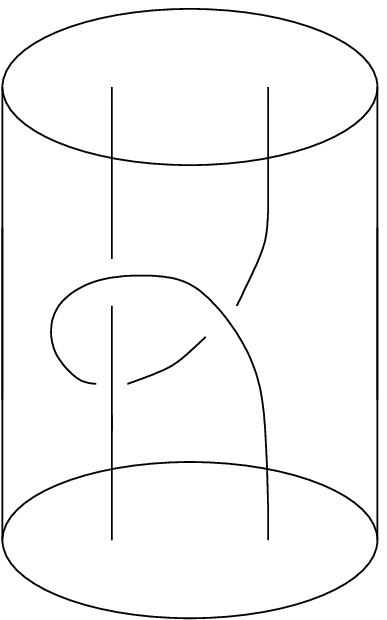, height=2cm}}
\hspace{1.5cm}
\subfigure[]{\epsfig{file=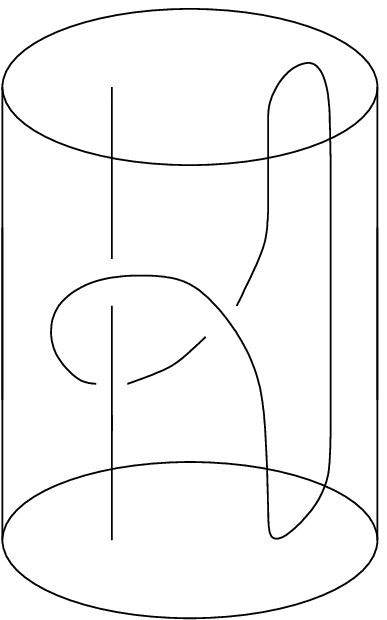, height=2cm}}
\caption{A string-link and its deformation closure.}
\label{stringlink}
\end{figure}

\smallskip
We now turn our attention to defining string-links in a homology ball $B^M$. These are the type of tangles we will be concerned with.

Fix a collection of points $p_1 < p_2 < p_3 < \cdots $ on the $x$-axis of $D^2$.  These induce  sets of points on $\partial ( \dti )$ which we call the \textit{standard  points}.
By a \textit{string-link of n components}  $\sigma \subset B^M$  we mean a tangle
$\sigma : \cup_{i=1}^{n} I_i \rightarrow B^M $ such that under the identification of the boundary with $\partial (\dti)$ we have  $\sigma |_{I_i} (j) = p_i \times j$, for $j=0,1$, where $p_i$ is the $i$-th standard point and $I$ is the closed unit interval. Note that a string-link has no closed components.

The two sets of {\em standard points} of a string link are the subsets of standard points on $D^2 \times \{ j \}$, $j=0,1$ which are the boundary points of a string-link component.

Given a string-link $\sigma \subset \dti$, we can change a given set $A$ of (interval) components of $\sigma$ into  $S^1$ components by constructing non-intersecting paths on $\pa{\dti}$  between the two endpoints of each  interval component, and pushing these paths and the endpoints of the components of $A$ slightly into the interior of $\dti$.  This gives a tangle which we call  the  {\em deformation closure}  of $\sigma$ with respect to $A$.  An example is given in Figure~\ref{stringlink}, where the deformation closure is with respect to the right hand component.

We can {\em represent} any string-link $\sigma^{\prime} \subset B^M$ by a string-link $\sigma \subset \dti$ with a specified set of framed surgery components, where $\sigma^{\prime}$ is obtained by carrying out surgery on the deformation closure of the surgery components of $\sigma$.

\smallskip

A {\em parenthesization}  of a set of standard points is a bracketing of that set (for example $(( p_1 p_2 )( p_3 ( p _4 p_5 )))$.  We call the parenthesization $(((( p_1 p_2 ) p_3 ) \cdots ) p_n )$ the {\em canonical parenthesization}.
A {\em parenthesization}  of a string-link  is a parenthesization of its two sets of standard points.

Note that a parenthesization on a string-link in $\dti$ induces one on any string-link in $B^M$ it represents.

\begin{Def}
A {\em manifold string-link} is a  canonically  parenthesized, \linebreak coloured, framed string-link in $\dti$ with a set of linking components $X_L$ and surgery components $X_M$.
\end{Def}

We say that a manifold string-link is {\em regular} if the linking matrix of its surgery components is invertible (so  surgery yields a string-link in  a rational homology ball).

\smallskip

We will now define some actions on the set of string-links which we will make use of later.  As these are well known and somewhat fiddly to define, we gloss over the technical details and rely upon the reader's intuition.

 If two string-links $\sigma_1$ and $\sigma_2$ in $\dti$ have the same number of components then we may form a product $\sigma_1 \cdot \sigma_2$ in the usual way by ``putting $\sigma_2$ on top of $\sigma_1$''.
 If the string links are parenthesized or coloured we require that the parenthesization or colourings match on the two disks  identified under the composition.

We also define $\sigma_1 \otimes \sigma_2$ to be the string link obtained by ``placing $\sigma_2$ to the right of $\sigma_1$''.

Let $T \subset B^M$ be an $X$-coloured tangle and let $A \subset X$.
 Define $\varepsilon_A (T)$ to be the tangle obtained from $T$ by deleting all of the components with colours in $A$.
Further let $B$ be a set disjoint from $X$  and let $S \subset X \times B$ such that any element of $X$ or $B$ appears in at most one pair $(x , b) \in S$.
We define $D_{S}(T)$ to be the coloured tangle obtained from $T$ by, for each $(x,b) \in S$, doubling the $x$-coloured component and colouring the double with $b$.
When dealing with string-links we may have to isotope them so that the endpoints lie on the appropriate standard points.

\section{Milnor's $\mu$-invariants}
Recall that given a ring $R$, a {\em R-homology sphere} is a 3-manifold $M$ such that $H_{q}(M;R) =H_{q}(S^3 ;R)$, for all integers $q$.  Similarly a {\em R-homology ball } is a 3-manifold $B^M$ with boundary $\partial (B^M ) = S^2$ such that $H_{q}(B^M;R) =H_{q}(B^3 ;R)$, for all  $q$, where $B^3$ is the 3-ball.
If  $R= \mathbb{Z}$  we do not specify the ring and just write {\em homology sphere} or {\em homology ball}.

Let $\sigma$ be a $l$-component string-link in a homology ball $\bsig$ with a fixed identification of $\partial (\bsig)$ with $\partial (\dti )$ (ie. $\bsig$ is a homology cylinder over $D^2$) and let  $N(P)$ be a regular neighbourhood of the set of standard points and  $N( \sigma )$  a regular neighbourhood of the string-link. Then there are two inclusion maps
 $i_j :D^2 - N(P) \hookrightarrow \bsig  - N( \sigma )$ for $j=0,1$,
 where the map $i_j$ sends $D^2 - N(P)$ to the image of $D^2 \times \{ j \} - N(P)$ under the identification of  $\partial (\bsig)$ with $\partial (\dti )$.
We use the $i_j$ to induce certain isomorphisms as follows.

Let  $G$ be any group. The \textit{lower central series}, $G_q$ is defined inductively by $G=G_1$ and  $G_{q+1} = [G,G_q]$.

\begin{StThe}[\cite{St:65}]
Let $\map{h}{A}{B}$ be a homomorphism of groups, inducing an isomorphism $\h{1}{A} \cong \h{1}{B} $ and an epimorphism from $\h{2}{A}$ onto $\h{2}{B}$.  Then, for finite $q$, $h$ induces an isomorphism $ \quoq{A} \cong \quoq{B}$.
\end{StThe}

A Mayer-Vietoris calculation and a standard  application of Stallings' theorem gives the following result.

\begin{Prop} \label{prop:isoms}
$(i_j)_*$, $j=0,1$, induces isomorphisms
\[ \fquo{\fun{(D^2 \times \{j\})- N(P)}} \cong \fquo{\fun{\bsig - N(\sigma)}}. \]
\end{Prop}

Let $F(l)$ be the free group on generators $x_1, \ldots , x_l$.  We will also denote the image of $x_i$ in the quotient group $\quoq{F(l)}$  by $x_i$ and the induced maps on the lower central series coming from proposition~\ref{prop:isoms}  by $(i_j)_*$, $j=0,1$.
Since  we can identify $\fun{(D^2 \times \{j\}) - N(P)}$ with $F(l)$, we have isomorphisms
\[ \fquo{F(l)} \overset{(i_0)_*}{\longrightarrow} \fquo{\fun{\bsig - \sigma}} \overset{(i_1)_*}{\longleftarrow} \fquo{F(l)}, \]
and the composition $(i_1)_{*}^{-1} (i_0)_*$ gives a map
$ SL(l) \rightarrow Aut(\quoq{F(l)})$, where  $SL(l)$ is the set of string-links of $l$ components in a given homology ball $\bsig$.
It is not difficult to see that we in fact get a map
\[ \map{\art_q}{SL(l)}{Aut_0 (\quoq{F(l)})}, \]
where $Aut_0 (\quoq{F(l)})$ is the subgroup of $Aut (\quoq{F(l)})$ consisting of all automorphisms which map $x_i$ to a conjugate of itself and leaves the product $x_1 x_2 \cdots x_l$ of the generators fixed.
We call the map $\art_q$ the $q$-th \textit{Artin representation}.

The $i$-th \textit{longitude} $\lambda_i \in \quoq{F(l)}$ of a string-link $\sigma$ is defined in the following way.
Take a double of the $i$-th component of the string-link. This determines an element in the fundamental group of the complement,  under $(i_1)_*^{-1}$ this gives an element in $\quoq{F(l)}$ which we call the longitude.
We have
\[ \art_q (\sigma)(x_i)=\lambda_i x_i \lambda_{i}^{-1}, \]
where $\lambda_i \in \quoq{F(l)}$ is the $i$-th longitude of $\sigma$.

Note that our longitudes are determined by the (black-board) framing and are not necessarily null-homologous.  It is easy to modify the content of this chapter should we insist that the longitudes are null-homologous, or we could just use the zero framing.

\begin{Def} We say that a string-link $\sigma$ has \textit{Milnor filtration n}, if all its longitudes are trivial in $\frac{F(l)}{F(l)_n}$.
\end{Def}

\smallskip

Let $\pow{l}$ be the ring of formal power series in non-commuting variables $\ser{X}{l}$. The {\em Magnus expansion} is the homomorphism
\[ \map{\mu}{F(l)}{\pow{l}} \]
defined on the generators of the free group by $\mu (x_i) = 1 + X_i$.

\begin{Def}
The {\em $\mu$-invariants}  of a string-link $\sigma$ in an integral homology ball are the coefficients of the monomials in the $X_i$ of the Magnus expansion of the $i$-th longitude $\lambda_i \in \fquo{F(l)}$.  Explicitly, the $\mu$-invariant of {\em length} $n+1$ of $\lambda_i$ is
\[
\mu_{j_1 ,j_2 ,   \ldots , j_n ;i} = \text{Coeff}( X_{j_1} X_{j_2} \cdots X_{j_n} , \mu (\lambda_i))
\]
where $n \leq q$ and $\lambda_i \in \fquo{F(l)}$.
\end{Def}

It is well known that the longitudes $\lambda_i$ of $\sigma$ are trivial in $\frac{F(l)}{F(l)_n}$, that is $\lambda_i$ is of Milnor filtration $n$, if and only if all $\mu$-invariants of length $\leq n$ vanish.

\section{The Algebras} \label{sec:algebras}
The algebras we need are amalgamations of the usual algebras $\mathcal{A}$ and $\mathcal{B}$ from the theory of finite-type invariants (see~\cite{BN:95:2, BN:95}).

\begin{Def}
Let $X,Y$ be finite disjoint sets.
 Then $\algyz$ is the space of formal $\mathbb{Q}$-linear combinations of uni-trivalent graphs whose trivalent vertices are oriented and whose univalent vertices are either coloured by elements of a set $Y$ or lie on the oriented coloured 1-manifold $(\cup_{x \in X} I_x )$,  which is called the \textit{skeleton},  modulo the STU, IHX and AS relations shown in figure~\ref{fig:relations}.
\end{Def}

\begin{figure}
\[
\begin{array}{c}
\mathrm{STU:} \hspace{0.5cm} \STU{1cm}{0.5cm} \hspace{1cm} \mathrm{AS:} \hspace{0.5cm} \AS{1cm}{0.5cm} \\
\vspace{0.3cm} \\
\mathrm{IHX:} \hspace{0.5cm} \IHX{1cm}{0.5cm}
\end{array}
\]
\caption{The STU, AS and IHX relations.}
\label{fig:relations}
\end{figure}

Note that we allow trivalent graphs and the possibility that $Y=\emptyset$.

We denote the subspace of $\algyz$ such that every connected component has a univalent vertex and all univalent vertices lie on the skeleton by $\mathcal{A}(\uparrow_X)$ and the subspace $\algez$ by $\mathcal{B}(Y)$.

The \textit{degree} of a uni-trivalent diagram is half of its number of vertices
and we say  an element of $\algyz$ is {\em connected} if it is a $\mathbb{Q}$-linear combination of connected uni-trivalent graphs.

\smallskip

\begin{figure}
\begin{center}
\epsfig{file=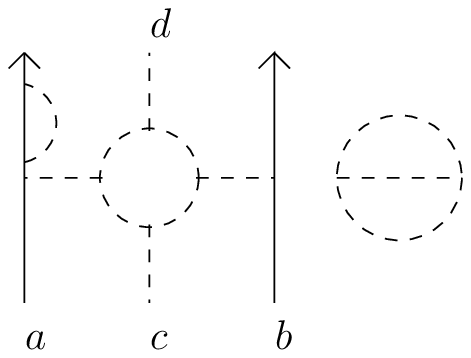, height=2cm}
\raisebox{10mm}{$\cdot$}
\epsfig{file=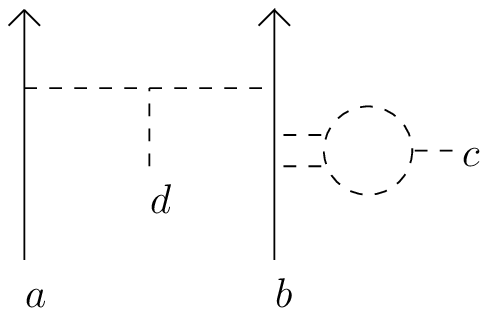, height=2cm}
\raisebox{10mm}{$=$}
\epsfig{file=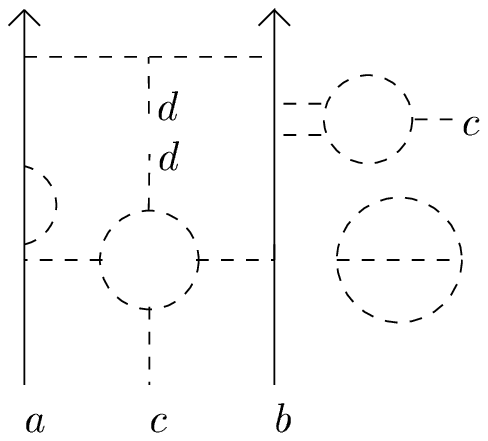, height=2cm}
\end{center}
\caption{An example of multiplication.}
\label{exmult}
\end{figure}

Let $D_1, D_2 \in \algyz$ then there is a  {\em product} $D_1 \cdot D_2$ given by the linear extension of the process of stacking the skeleton of $D_1$ on top of $D_2$ in such a way that the colours of the two skeletons match and taking the disjoint union of any trivalent components.
An example of the multiplication is given in figure~\ref{exmult}.

There is also a notion of a coproduct $\cp$ in $\algyz$ which is the obvious extension of the usual coproduct of $\mathcal{A}$ (see~\cite{BN:95:2}). In fact this makes  $\algyz$ into a graded co-commutative Hopf algebra where the grading is by the degree.
We denote the degree $n$ part by $\algyzn $ and, by abuse of notation, its graded completion again by $\algyz$.
 The primitives (ie. the elements such that $\cp (D) = 1 \otimes D + D \otimes 1$) of the algebra are the connected elements.

\smallskip

We will now look at some maps between these algebras. All of these properties hold since they hold in $\mathcal{A}$ and $\mathcal{B}$.

\begin{figure}
\[ \chiso_{\{b\}}
\left( \begin{array}{c}
 \epsfig{file=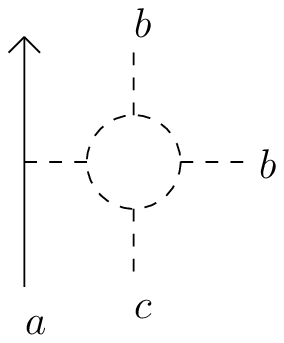, height=1.8cm}
\end{array} \right)
= \frac{1}{2}
\left( \begin{array}{ccc}
 \epsfig{file=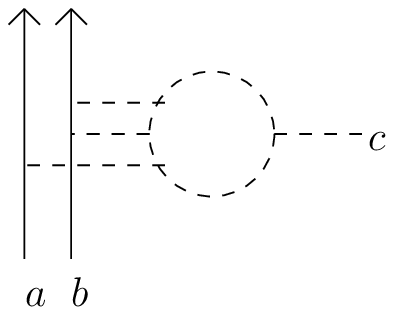, height=1.8cm}
 & \raisebox{0.9cm}{+}
  & \epsfig{file=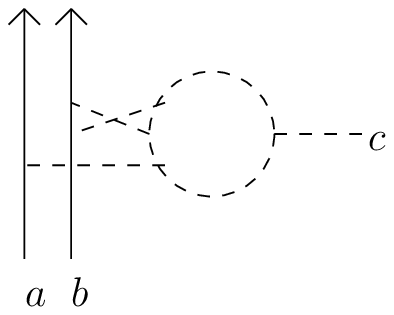, height=1.8cm}
\end{array} \right)
\]
 \caption{An example of the map $\chi_{ \{ b \} }$. }
\label{exchi}
\end{figure}

Let $Y^{\prime} \subset Y$. Define a map
\[
\chiso_{ Y^{\prime} } : \algyz \rightarrow \algchiimage
 \]
by the linear extension of the process of adding $Y^{\prime}$ coloured skeleton components and
taking the average of all ways of placing the $ Y^{\prime}$ labeled univalent vertices on $\uparrow_{ Y^{\prime}}$.
See figure~\ref{exchi} for an example of this map.

In fact $\chiso$ descends to a coalgebra isomorphism and we denote its inverse by $\sigma$.

If $X^{\prime} \subset X$, $Y^{\prime} \subset Y$  and $A=X^{\prime} \cup Y^{\prime}$. The map
$\varepsilon_{A} : \algyz \rightarrow \mathcal{A}( \uparrow_{X -X^{\prime}}, Y -Y^{\prime})$
is defined by setting every uni-trivalent graph with a uni-valent vertex on a $X^{\prime}$ coloured skeleton component or  with a $Y^{\prime}$ coloured vertex equal to zero.

\begin{figure}
\[ D_{\{(a,b) , (c,d) \}}
\left( \begin{array}{c}
 \epsfig{file=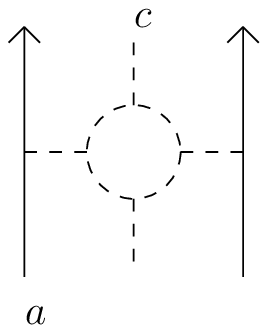, height=1.8cm}
\end{array} \right)
=
\begin{array}{c}
 \epsfig{file=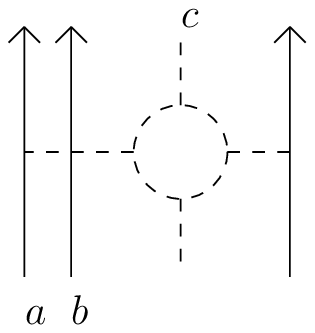, height=1.8cm}
 \raisebox{0.9cm}{+}
 \epsfig{file=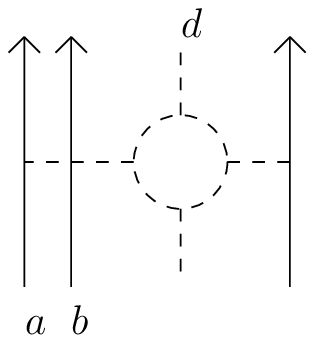, height=1.8cm}
 \raisebox{0.9cm}{+}
 \epsfig{file=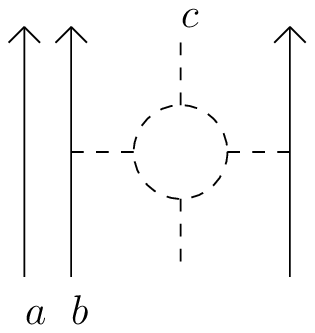, height=1.8cm}
 \raisebox{0.9cm}{+}
 \epsfig{file=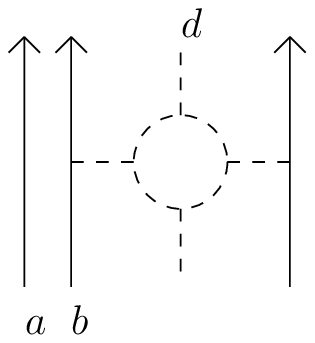, height=1.8cm}
\end{array}
\]
 \caption{An example of the map $D_{S}$.}
\label{exd}
\end{figure}

Let $B$ be some set  disjoint from both $X$ and $Y$ and let
$S \subset (X \cup Y) \times B$ such that any element of $X,Y,B$ appears in at most one pair $(a , b) \in S$. Then define
$D_{S} $ to be the linear extension of the operation which to each element $(a,b) \in  S$
 either, if $a$ is a label of a skeleton component, gives the sum of all ways of lifting the vertices lying on the $a$-coloured component  over    the component and its $b$-coloured double   and, if $a$ is the colour of a univalent vertex, is the sum of all ways of substituting the colour $a$ by $b$  (see figure~\ref{exd} for an example).

\medskip

We will be interested in two particular quotients of $\algyz$ which were  defined in \cite{BN:95} and \cite{HM:00} for the algebra $\mathcal{A}$.

Define $\algyzt$ to be the quotient of $\algyz$ by the ideal generated by all relations which set non-simply connected uni-trivalent graphs  equal to zero.  The connected elements are called  {\em trees} .
We will denote the connected part (ie. the primitives) of $\mathcal{B}^t (Y)$ by $\mathcal{C}^t(Y)$.

Also define $\algyzh$ to be the quotient of $\algyzt$ by the ideal generated by all relations which set connected uni-trivalent graphs with more than one univalent vertex either lying on the same skeleton component  or being labeled by the same colour, equal to zero.  We call this quotient the {\em homotopy quotient}.

It follows from \cite{BN:95}  that $\algyzh$ is a quotient of $\algyzt$ and $\chiso$ descends to  isomorphism on  these algebras.

\smallskip

There is a well known map (see eg. \cite{GL:P1, HM:00}) which relates trees to Lie algebras.
Let $ \text{Lie} (l) = \oplus_{n \geq 1} \text{Lie}_n (l)$, be the free $\mathbb{Q}$ Lie algebra on $l$ generators $\ser{X}{l}$.
Also let $\mathcal{C}^t(Y,a)$ be the subspace of $\mathcal{C}^t(Y \cup \{ a \})$ consisting of connected elements in which every uni-trivalent graph has exactly one univalent vertex coloured by some $a \notin Y$.

\begin{figure}
\[
\epsfig{file=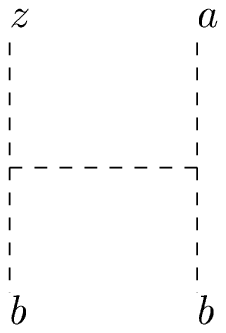, height=2cm}
\raisebox{0.9cm}{$\rightarrow$}
\epsfig{file=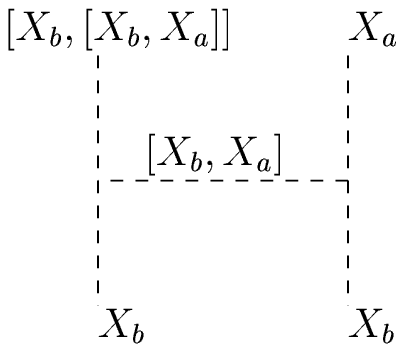, height=2cm}
\raisebox{0.9cm}{$\rightarrow [X_b ,[X_b , X_a]]$}
\]
\caption{An example of the isomorphism
$\mathcal{C}^t_3 (\{ a,b \},z) \rightarrow \text{Lie}_3 (2)$.} \label{fig:primlie}
\end{figure}

Fix a bijection between the colouring set $Y$ and the generators  $\ser{X}{l}$ of the free Lie algebra, where $|Y|=l$. Then given some element $D \in \mathcal{C}^t_n(Y,a)$  label the edges ending in a $Y$-coloured univalent vertex with the corresponding generator of the Lie algebra.
Now assign an element of the Lie algebra to each unlabelled edge according to the rule that whenever an unlabelled edge meets two edges labelled by $X$ and $X^{\prime}$ in $ \text{Lie} (l)$ (in the direction of the orientation) assign the commutator $[X,X^{\prime}]$ to that edge.  This labels the edge coloured by $a$ and we take this to be our element of $\text{Lie}_n (l)$.
See figure~\ref{fig:primlie} for an example.
It is not hard to see that this  gives an isomorphism from $\mathcal{C}^t_n(Y,a)$ to $\text{Lie}_n (l)$.

\begin{figure}
\[
\epsfig{file=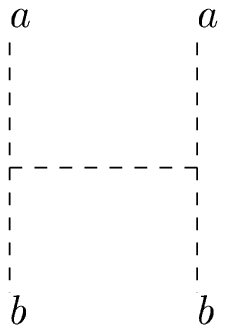, height=2cm}
\raisebox{0.9cm}{$\rightarrow$}
\hspace{2mm}
\epsfig{file=mufig/is1, height=2cm}
\raisebox{0.9cm}{+}
\hspace{2mm}
\epsfig{file=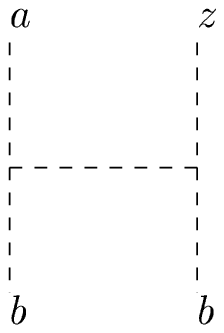, height=2cm}
\raisebox{0.9cm}{$\rightarrow$}
\hspace{2mm}
\raisebox{0.9cm}{2}
\hspace{1mm}
\epsfig{file=mufig/is1, height=2cm}
\raisebox{0.9cm}{$\rightarrow 2[X_b ,[X_b , X_a]]$}
\]
\caption{An example of
$j_z : \mathcal{C}^t_3 (\{ a,b \}) \rightarrow \text{Lie}_3 (2)$.
} \label{fig:mapj}
\end{figure}

Finally, we define a map
$ j_y : \mathcal{C}^t_n(Y) \rightarrow \text{Lie}_n (l)$ for $y \in Y$
by summing over all of the ways replacing exactly one of the $y$-coloured vertices by some $a \notin Y$ and then using the above map to get an element in $\text{Lie}_n (l)$.
An example is given in figure~\ref{fig:mapj}.

\bigskip

Before continuing we briefly review some relevant of the Kontsevich integral.

\begin{itemize}
\item For our purposes the Kontsevich integral $Z$ is an $\mathcal{A}(\uparrow_{X})$ valued universal finite-type invariant of $X$-coloured framed parenthesized tangles and the degree $n$ part $Z_n$ of $Z$   is a degree $n$ finite-type invariant.

\item Let $\pi^{h} : \mathcal{A}(\uparrow_{X}) \rightarrow \mathcal{A}^h(\uparrow_{X})$ be projection. Then by \cite{BN:95}, $\pi^h \circ Z$ is a well defined invariant of ($X$-coloured framed parenthesized) string-links  up to link-homotopy, where {\em link-homotopy} is an equivalence relation which allows ambient isotopy and each component of the tangle to pass through itself.

\item Let $T$ and $T^{\prime}$ be tangles then $Z(D_{S}(T)) = D_{S}(Z(T))$,
$Z(\varepsilon_{A}(T)) = \varepsilon_{A}(Z(T))$ and $Z(T \cdot T^{\prime}) = Z(T) \cdot Z(T^{\prime})$.

\item $Z(T) \in \mathcal{A}(\uparrow_{X})$ is group-like and so can be written as $\exp ( C )$ where $C \in \mathcal{C}(X)$ is connected.
\end{itemize}

\section{The \AA rhus Integral} \label{sec:aarhus}

The \AA rhus integral, $Z^M$,  was introduced by Bar-Natan, Garoufalidis, Rozansky and Thurston in the series of papers \cite{BGRT:AI, BGRT:AII, BGRT:AIII} as a universal finite type invariant of rational homology 3-spheres.
In this series it was remarked that it  extends to an invariant of  links in rational homology spheres.
In this section we define \AA rhus integral.
The reader is referred to the \AA rhus trilogy  for a thorough exposition  of the invariant.

\medskip

The {\em pre-normalized \AA rhus integral} of regular manifold string-links  $\PZM$ is defined by the following composition:
\[
\PZM : RMSL \overset{\z}{\longrightarrow} \Alm \overset{\sigma_{X_M}}{\longrightarrow} \alglm \overset{\fgi}{\longrightarrow} \algle
\]
where:
\begin{itemize}
\item  $RMSL$ is the set of regular manifolds string-links with linking components coloured by $X_L$ and  surgery components coloured by $X_M$.
\item $\z \overset{def}{=} \nu^{\otimes |X_L \cup X_M |} \cdot D_{\{-\}\times X_L \cup X_M} (\nu ) \cdot Z$, is the Kontsevich integral as normalized in \cite{LMMO:99}.
\item $\fgi$ is  {\em formal Gaussian integration} with respect to the variables $X_M$. It is described below.
\end{itemize}

\begin{Def}
The {\em \AA rhus integral} of a regular manifold string-link, $\sml$ is given by
\[
\ZM(\sml ) = \PZM (U_+ )^{- \sigma_{+}} \cdot \PZM (U_- )^{- \sigma_{-}} \cdot  \PZM(\sml )
\]
where $\sigma_{\pm}$ is the number of $\pm$ve eigenvalues of the linking matrix of $\varepsilon_{X_L} (\sml) $ and $U_{\pm}$ is the unknot with framing $\pm 1$.
\end{Def}

We will now go on to define formal Gaussian integration.
Let $D_1, D_2 \in \algrs$, define
\[
\bracy{D_1}{D_2} =
\left(
\begin{array}{l}
\text{sum of all ways of gluing all legs labeled } \\ y \text{ on } D_1 \text{ with  all of the legs labeled } y \text{ on } \\ D_2 \text{, for all the colours } y \in Y .
\end{array}
\right) ,
\]
 where this sum is non-zero only if the number of $y$-coloured legs of $D_1$ equals the number of $y$-coloured legs of $D_2$, for all $y \in Y$.

\begin{figure}
\[
\left\langle \begin{array}{c}
 \raisebox{0.7cm}{$\frac{1}{2}$}
 \hspace{1mm}
 \epsfig{file=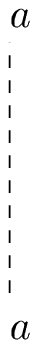, height=1.4cm}
\raisebox{0.7cm}{$+$}
\hspace{2mm}
\epsfig{file=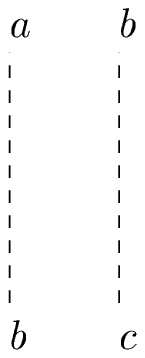, height=1.4cm}
\hspace{2mm}
,
\hspace{2mm}
\raisebox{0.7cm}{$4$}
\hspace{1mm}
 \epsfig{file=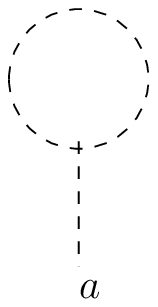, height=1.4cm}
\hspace{1mm}
\raisebox{0.7cm}{$+$}
\hspace{1mm}
 \epsfig{file=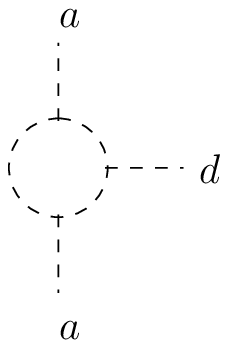, height=1.4cm}
 \raisebox{0.7cm}{$+$}
 \epsfig{file=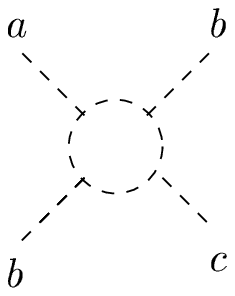, height=1.4cm}
\end{array} \right\rangle_{ \{ a,b,c \} }
\begin{array}{c}
\raisebox{0.7cm}{$=$}
 \hspace{3mm}
 \epsfig{file=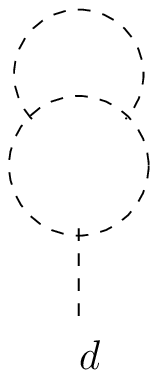, height=1.4cm}
 \hspace{1mm}
 \raisebox{0.7cm}{$+$}
 \hspace{2mm}
 \raisebox{0.7cm}{$2$}
 \hspace{1mm}
 \epsfig{file=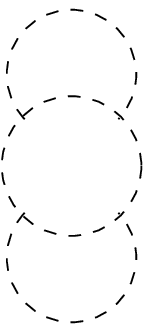, height=1.4cm}
\end{array}
\]
 \caption{An example of the pairing $\bracs{D_1}{D_2}$.}
\end{figure}

\medskip

It is a well known and easily seen fact that for an $X$-coloured tangle $T$, the Kontsevich integral $\z $ may be written in the form
\[ \sigma (\z (T)) = \exp_{\udot} (\sum_{x,y \in X} \frac{1}{2} l_{xy} \; x \frown y ) + (\text{other stuff}), \]
where $(l_{xy})$ is the linking matrix of $T$.  Recall that the degree 1 elements (which look like $\frown$) are called {\em struts}.

Therefore, given a regular manifold string-link $\sml$, with linking components $X_L$ and surgery components $X_M$, we can write
\[
\sigma_{X_M} \z (\sml) = \exp_{\udot} (\sum_{x,y \in X_M} \frac{1}{2} l_{xy} \; x \frown y) \cdot P \overset{def}{=} \exp_{\udot} (Q/2)  \cdot P,
\]
where now $(l_{xy})$ is the linking matrix of $\varepsilon_{X_L}(\sml)$.
Since $\sml$ is regular, $(l_{xy})$ is invertible and so we can define:
\[
Q^{-1} = \sum_{x,y \in X_M} l^{xy} \; x \frown y
\]
where $(l^{xy})$ is the inverse matrix of $(l_{xy})$.

Writing $\sigma_{X_M} \z (\sml) \overset{def}{=}  \exp_{\udot} (Q/2)  \cdot P$, we define {\em formal Gaussian integration} as:
\[
\fgi  \exp_{\udot} (Q/2)  \cdot P = \bracxm{\exp_{\udot} (-Q^{-1} /2)}{P}.
\]

\medskip

It  is known (\cite{BGRT:AII}) that $\PZM$ is invariant under isotopy and a handle slide of any component around a surgery component, and  $\ZM$ is invariant under stabilization on the surgery components and so $\ZM$ is an invariant  of string-links in a rational homology ball.
Summarizing this we have:

\begin{Prop}
The \AA rhus integral $Z^M$ is an invariant of framed parenthesized string-links in rational homology balls.
\end{Prop}

\medskip

At this point we fix some notation.
Let $\pi^t : \algxy \rightarrow \algxyt$ be projection.  Then
$\ZMT \overset{def}{=} \pi^t \circ \ZM$,
$\bracy{-}{-}^t \overset{def}{=} \pi^t \circ \bracy{-}{-}$
and so on.
We use similar notation for the projection $ \pi^h : \algxy \rightarrow \algxyh$

\section{The \AA rhus Integral and the $\mu$-invariants}

Let $\sml$ be a manifold string-link with the canonical parenthesization such that  the determinant of the linking matrix of the surgery components is $\pm 1$ (so $\sml$ represents a string-link in an integral homology ball).
Further, for convenience, we set $X_L = \{1, \ldots , l \}$ ,  $X_M = \{l+1, \ldots , l+m \}$ and assume that the components of the manifold string-link have numerically increasing colours from left to right.
We call such a manifold string-link {\em nice}.

At times we will need to add an extra  linking component to the manifold string-link.  We will add this component to the left of the others and colour it with $0$.
We denote the new colouring set $X_L \cup \{ 0 \}$ by $X_L +1$.

The extra $0$-coloured component is going to correspond to a longitude of the string-link and as such is only considered up to link homotopy.  Consequently, rather than working with the algebra $\algtl1e$, we add an additional homotopy relation on the colour $0$, and we call the resulting algebra $\alghtl1e$.

Given a set of colours $X$, let  $1_X $, be the trivial tangle of $|X|$ components coloured by $X$. When $X$ contains only one element, $x$ say,  we will just write $1_x$.

In this section we consider the longitudes as elements
$\lambda_i = \lim_{ \hspace{-0.6cm}\raisebox{-1mm}{$\longleftarrow$}} \lambda_i^{(n)}$ of the nilpotent completion
$ \widehat{F(l)} = \lim_{ \hspace{-0.6cm}\raisebox{-1mm}{$\longleftarrow$}} \quo{F(l)}$,
where $\lambda_i^{(n)} \in \quo{F(l)}$.

The reader is referred to \cite{HM:00} for the motivation behind the formula in the following lemma.

\begin{Prop} \label{mu:lambda}
Let $\sml$ be a nice  manifold string-link   and let  $\lambda_i$, $1 \leq i \leq l$, be its $i-th$ longitude regarded as a pure braid. Then
\begin{equation} \label{mu:form}
Z^{M;h,t}(\lambda_i \otimes 1_{X_M}) = \pi^{h,t} (Z^M (1_0 \otimes \sml )^{-1} (D_i Z^M ( \sml ))^{b_i})
\end{equation}
where $b_i = Z^M (\beta_i)$ and  $\beta_i$ is the braid coloured by $\{0, \ldots , l+m \}$ inducing the permutation $(i-1 \; i-2 \cdots  1 \; 0)$, $a^b$ denotes the conjugation  $bab^{-1}$, $D_i = D_{\{(i,0)\}}$ and $\pi^{h,t}$ is projection onto $\alghtl1e$.
\end{Prop}

\begin{Remark}
In formula~\ref{mu:form} we are assuming that $\sml$, $\lambda_i \otimes 1_{X_M}$ and $1_0 \otimes \sml$ have the canonical parenthesization and $\beta_i$ has the canonical parenthesization on the bottom and the `$i$-th double of the canonical parenthesization' on the top.
\end{Remark}

\begin{Remark}
 Since $\Zmht_0 = \Zmht$ we need only consider the pre-normalized \AA rhus integral. Also note that $\z^t =Z^t$.
\end{Remark}


\medskip

We need a few technical lemmas to prove the proposition.

\begin{Lem} \label{mu:double}
Let $\sml$ be a manifold string-link, $i \in X_L$ and $D_i = D_{\{(i,0)\}}$. Then
\[Z^{M}(D_i (\sml )) = D_i(Z^{M}(\sml )).\]
\end{Lem}

\begin{proof}
\[
\begin{split}
\pzm{D_i ( \sml )} &=  \sideset{}{_{X_M}^{FG}}\int  \sigma_{X_M} \z (D_i ( \sml)) \\
                &=  \sideset{}{_{X_M}^{FG}}\int  D_i( \sigma_{X_M}( \z( \sml ))) \\
                &=  D_i \left( \sideset{}{_{X_M}^{FG}}\int  \sigma_{X_M} \z (\sml ) \right) \\
                &= D_i(\pzm{ \sml})
\end{split}
\]
where the second equality is a standard property of the Kontsevich integral.
The third follows since $i \in X_L$ and the formal Gaussian integration is with respect to the variables $X_M$.

The result follows since $D_i$ respects multiplication.
\end{proof}

\begin{Lem} \label{mu:glue}
Let $P \in \algyz$ be of degree $n$ and contain no struts both of whose univalent vertices are coloured by elements of $Y$, and let $Q \in \algyz$ consist entirely of struts coloured by $Y$.  Then if $\bra{Q}{P}{Y}^t$ is non-zero, it is of degree at least $ n-[ \frac{n}{2} ]$, where  $[a]$ is the integer part of $a$.
\end{Lem}

\begin{proof}
Assume that  $\bra{Q}{P}{Y}^t$ is non-zero.
We need to find the minimum possible degree of $\bra{Q}{P}{Y}^t$, where $P$ and $Q$ vary over all suitable elements of $\algyz$.

Suppose we are given elements $P$ and $Q$ so that the degree of $\bra{Q}{P}{Y}^t$ is minimal.
We can assume that $P$ is simply connected.
If there are any trivalent vertices in $P$ then we remove them by identifying two of the edges incident to the trivalent vertex giving a single edge and adding a $X$-coloured univalent vertex at the end of the third edge (note we use the hypothesis that $P$ is a tree). As this does not change the degree of $P$ or $\bra{Q}{P}{Y}^t$, we see we can assume that $P$ consists entirely of struts.

Given such a $P$ and remembering that the struts in $P$ have at most one $Y$-coloured vertex, it is easy to see that the maximum number of $Y$-coloured struts which may be glued in is $[ \frac{n}{2} ]$, giving the result.
\end{proof}

Note that if $P$ and $Q$ are as in the above lemma and if $\bra{Q}{P}{Y}$ is non-zero then it is of degree at least $ n-[ \frac{3n}{4} ]$.





\begin{figure}
\begin{center}
 \epsfig{file=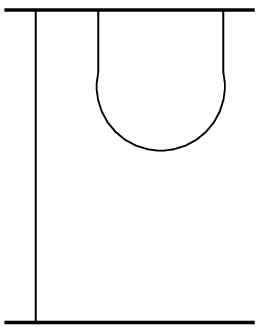, height=1.8cm}
\hspace{1cm}
\epsfig{file=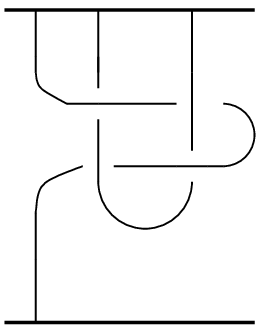, height=1.8cm}
\caption{Two tangles which differ by a pure braid.}
\label{mu:differ}
\end{center}
\end{figure}

\begin{Def}
We say that two tangles  $T$ and $T'$ {\em differ by a pure braid $p \in PB_{n+1}$} if $T'$ can be obtained from $T$ by replacing a copy of $D^2 \times I$ which intersects $T$ in a trivial string-link with the pure braid $p$  (see figure~\ref{mu:differ}).
\end{Def}

We need the following result of Stanford.
\begin{The}[Stanford \cite{St:96}]
Let  $T$ and $T'$ be two tangle which differ by a pure braid $p \in PB_{n+1}$.  Then for any finite type invariant, $v$, of degree less than $n+1$ we have $v(T)=v(T')$.
\end{The}

\begin{Lem} \label{mu:uni}
Suppose Y is the disjoint union of compact 1-manifolds $Y_0$, $Y_L$, $Y_M$ where $Y_M$ consists entirely of copies of $S^1$.  Further suppose that $T_i$,  $i=1,2$ , are two tangles which agree on $Y_L \cup Y_M$ and on each component of $Y_0$ the maps differ by an element in the lower central series $ \pi_1 (M-T_i |_{Y_L} )_{n+1}$ , where $M$ is the manifold obtained by surgery on $T_i |_{Y_M}$.  Then the images of $\pzm{(T_i)}$ in $\algspecial$, where the homotopy filtration is on the $Y_0$ components, agree.
\end{Lem}

\begin{proof}
First note that $\pi_1 (\dti - T_i |_{Y_L \cup Y_M}  )$ is generated by the meridians of $T_i |_{Y_L \cup Y_M}$.
This means that  there exists a ball which intersects $T_i |_{Y_L \cup Y_M}$ in a trivial string-link such that the map of the fundamental groups induced by the inclusion of the trivial string-link into $T_i |_{Y_L \cup Y_M}$ is surjective.

Now since $ \pi_1 (M-T_i |_{Y_L} )$ is also generated by the meridians of $T_i |_{Y_L \cup Y_M}$ and $T_1$ and $T_2$ differ by elements in $ \pi_1 (M-T_i |_{Y_L} )_{n+1}$, we see that $T_2$ can be obtained from $T_1$ by handle sliding around the $Y_M$ components, modifying the $Y_0$ components by homotopy and modifying $T_1$ (inside the ball described above) by pure braids in
$F(|Y_L \cup Y_M|)_{n+1} \subset F(|Y_L \cup Y_M|) \subset \text{PB}(|Y_L \cup Y_M|+1) $.

Stanford's theorem tells us  that the modification by the pure braids does not affect finite-type invariants of degree less than $n+1$ and since the homotopy relations are applied to the $Y_0$ components so tangles differing under these two moves have the same image under $\z_{\leq n}^{h,t}$.
Finally, formal Gaussian integration takes care of the handle slides and the result then follows  by lemma~\ref{mu:glue}.
\end{proof}

\begin{proof}[Proof of proposition~\ref{mu:lambda}]
Let $\lambda_{i}^{(n)}$ be a representative of the longitude $\lambda_i$ in $F(l)/F(l)_{n+1}$ which we regard as a pure braid of $l+1$ components.
Now the two tangles
$(1_0 \otimes \sml)(\lambda_{i}^{(n)} \otimes 1_M)$ and $\beta_i D_i(\sml) \beta_i^{-1}$
both represent the union of the manifold string-link, $\sml$, and the longitude and therefore, by Stallings' theorem, satisfy the conditions of lemma~\ref{mu:uni}.
Then
\[
Z^{M;h,t}_{< n -[\frac{n}{2} ]}((\lambda_{i}^{(n)} \otimes 1_M) \cdot (1_0 \otimes \sml ) )=
Z^{M;h,t}_{< n -[\frac{n}{2} ]}(\beta_i \cdot D_i (\sml ) \cdot \beta_{i}^{-1}).
\]

Since $\lambda_i^{(n)} \in F(l)/F(l)_{n+1}$
it only shares crossings with the linking components of $\sml$, therefore
 \[Z^{M;h,t}_{< n -[\frac{n}{2} ]}((\lambda_{i}^{(n)} \otimes 1_M) \cdot (1_0 \otimes \sml ) ) = Z^{M;h,t}_{< n -[\frac{n}{2} ]}(\lambda_{i}^{(n)} \otimes 1_M) Z^{M;h,t}_{< n -[\frac{n}{2} ]} (1_0 \otimes \sml ).\]
Similarly,
 \[Z^{M;h,t}_{< n -[\frac{n}{2} ]}(\beta_i \cdot D_i (\sml ) \cdot \beta_{i}^{-1}) =
 Z^{M;h,t}_{< n -[\frac{n}{2} ]}(\beta_i )Z^{M;h,t}_{< n -[\frac{n}{2} ]}(D_i (\sml ) )Z^{M;h,t}_{< n -[\frac{n}{2} ]}(\beta_{i}^{-1} ).\]

Finally solving for  $Z^{M;h,t}_{< n -[\frac{n}{2} ]}((\lambda_{i}^{(n)} \otimes 1_M)$, and letting $n$ tend to infinity gives the result.
\end{proof}

\medskip

Having found a formula for the \AA rhus integral of the longitudes we turn our attention to finding a formula for the Magnus expansion of the longitudes.

\begin{Def}
An {\em expansion} is a homomorphism  $J:F(l) \rightarrow \mathcal{P}(l)$ such that
$J(x_i) = 1 +X_i + (\text{higher order terms})$
, where $F(l)$ is the free group on the generators $x_1, \ldots , x_l$ and  $\mathcal{P}(l)$ is the ring of formal power series in  non-commuting variables $X_1, \ldots , X_l$.
\end{Def}

Clearly the Magnus expansion is an expansion in this sense.
We show that the left hand side of formula~\ref{mu:form} can be regarded as an expansion and then we apply the following result of Lin to write the Magnus expansion of the longitudes in terms of the \AA rhus integral.

\begin{Lem}[Lin \cite{Li:97}] \label{lem:lin}
Let $J$ be any expansion and $\mu$ be the Magnus expansion.  Then there exist a unique unipotent automorphism $\Psi : \mathcal{P}(l)\rightarrow \mathcal{P}(l)$ such that $\mu = \Psi \circ J$.
\end{Lem}

Recall that a map $\Psi$ is said to be {\em unipotent} if for all $a \in  \mathcal{P}(l)$ of degree $n$, $\Psi(a)=a + O(n+1)$.

\medskip

 $\Ahtlo$ is a graded co-commutative Hopf algebra whose space of primitives is isomorphic to $\Cht{X_L +1}$, the space  of connected elements of $\Bhtlo$.
 Let $X_i \in \Ahtlo$ denote the element of degree 1 which has a single chord between the skeleton components coloured by 0 and $i$.

Then $X_1 , \ldots , X_l$  generate a free non-commutative power series ring
$\mathcal{P}(l) = \mathcal{P}(X_1 , \ldots , X_l) \subset \Ahtlo, $
(as the primitives of $\Ahtlo$ are isomorphic to
$\Cht{X_L +1}$ and this is naturally decomposed as
 $ \Ct{X_L,0} \oplus \Ct{X_L} = \text{Lie}(l) \oplus \Ct{X_L}$
and the first summand corresponds to $\mathcal{P}(X_L)$ by the isomorphism described in section~\ref{sec:algebras}).

\begin{figure}
\begin{center}
 \epsfig{file=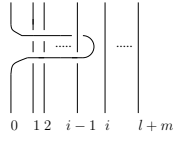, height=1.8cm}
\end{center}
\caption{The pure braid $ \sigma_{0,i} \otimes 1_M$.}  \label{mu:inj}
\end{figure}

Let  $SL(X_M )$ denote the monoid of string-links in $\dti$ which are coloured by $X_M$ and $PB(l+1)$ be the pure braid group on $l+1$ generators.
There is map $\iota: F(l) \rightarrow PB(l +1) \otimes SL(X_M)$ defined by the formula
$x_i \mapsto \sigma_{0,i} \otimes 1_M$
where $\sigma_{0,i}$ is the generator of the pure braid group which wraps the  $0$-th strand once around the $i$-th as in figure~\ref{mu:inj}.
The composition of this  with $ \Zmht $ gives a map $J:F(l) \rightarrow \alghtl1e$.

\begin{Lem}
The map $J :F(l) \rightarrow \alghtl1e  $ defined above is an expansion.
\end{Lem}

\begin{proof}
Let $x_i \otimes 1_M$ denote the generators of $F(l) \otimes id \subset PB(l+1) \otimes SL(X_M)$.
Then
\[
\zmht{x_i} = i \circ \zht{p_1(x_i)} = \Exp{Y_i}
\]
where $Y_i = X_i + O(2) \in \Ahtl1m$,
$i:\Ahtlo \rightarrow \Ahtl1m$ is  inclusion and $p_1: PB(l+1) \otimes SL(X_M) \rightarrow PB(l+1)$ is  projection onto the first component.

To prove the lemma we have to show that the image of $J$ lies in $\mathcal{P}(l)$.
But this follows since every diagram in $\zmht{x}$, where $x \in F(l) \otimes id$, must have a vertex lying on the skeleton component coloured $0$ since the removal of the $0$-coloured component trivializes the braid.
Thus $Y_i \in \Ct{X_L ,0} = \text{Lie}(l)$.
\end{proof}

Now applying Lemma~\ref{lem:lin} to Proposition~\ref{mu:lambda} gives:
\begin{The} \label{mu:mu}
Let $\sml$ be a nice manifold string-link.  Then
\[
\mu (\lambda_i) = \Psi \circ \pi^{h,t} (Z^M (1_0 \otimes \sml )^{-1} (D_i Z^M ( \sml ))^{b_i})
\]
where $b_i = Z^M (\beta_i)$ and  $\beta_i$ is the braid coloured by $\{0, \ldots , l+m \}$ inducing the permutation $(i-1 \; i-2 \cdots  1 \; 0)$, $a^b$ denotes the conjugation  $bab^{-1}$, $D_i = D_{\{(i,0)\}}$, $\pi^{h,t}$ is projection onto $\alghtl1e$ and $\Psi$ is a unipotent automorphism.
\end{The}

\begin{Remark}
Habegger and Masbaum's theorem in \cite{HM:00} relating the $\mu$-invariants of string-links in $\dti$ to the Kontsevich integral (which is obviously contained in theorem~\ref{mu:mu}) holds in the algebra $\Ahoxl1$, which we define to be $\Axl1$ with the homotopy relation applied to the colour 0, so there are no tree relations on the $X_L$-coloured components (this is true since \cite{HM:00}'s~lemma~12.5 only requires the homotopy relation and in their lemma~12.6, the homotopy relation ensures that the appropriate elements are trees).
However it is interesting to note that attempts by the author to remove the ``$t$'' in this section failed as the normalized Kontsevich integral $\z$ does not respect multiplication, necessitating the descent into $\alghtl1e$.
\end{Remark}

\section{The First Non-vanishing $\mu$-invariant}
In this section we give a generalization of the Habegger-Masbaum formula by expressing  the first non-vanishing Milnor invariants in terms of the first non-vanishing term of the tree part of the \AA rhus integral.
This result also appeared in \cite{Ha2000} with a different proof.

\medskip

Recall that $ \text{Lie} (l) $ is the free $\mathbb{Q}$ Lie algebra on $l$ generators $\ser{X}{l}$.
There is a canonical graded isomorphism of
$\oplus_{n\geq 1}(F(l)_n / F(l)_{n+1} ) \otimes \mathbb{Q}$ with $ \text{Lie} (l)= \oplus_{n \geq 1} \text{Lie}_n (l)$ (see \cite{MKS:book}).
Now if $\sigma$ has Milnor filtration $n$, we can consider the longitudes $\lambda_i$ as elements in $F(l)_n / F(l)_{n+1}$ and we denote the corresponding element in $\text{Lie}_n (l)$ by $\MU{i}{n}$.
We call the $\MU{i}{n}$  the {\em Milnor invariants of degree n}.

\begin{The} \label{mu:formulae}
Let $\sml$ be a nice manifold string-link representing a string-link $\sigma$. Then

(i) $\Zmht (\sml) = 1+ O(n)$ if and only if $\sigma$ is of Milnor filtration $n$,

(ii) the first non-vanishing Milnor invariants of the  string-link $\sigma$ determine and are determined by the first non-vanishing term of $\zmt{\sml} -1$ through the Habegger-Mausbaum formula:
\[ \mu_{i}^{(n)} (\sml) = j_i (\xi )  \]
where $\zmt{\sml} = 1 + \xi + O(n+1)$ and $j_i : \mathcal{C}^{t}_{n} (X_L) \rightarrow \text{Lie}_n (X_L) $ is the map described in section~\ref{sec:algebras}.
\end{The}

\begin{proof}

\noindent (i)
First suppose that $\zmt{\sml} = 1+O(n)$, then
$\pht (D_i(\zm{\sml})) = 1+O(n) $ giving $\pht ((D_i(\zm{\sml}))^{b_i}) = 1+O(n)$
(since the lower degree terms of the conjugating $b_i$'s cancel).

Also we have  $\pht (\zm{1_0 \otimes \sml}^{-1}) = 1 +O(n)$.

Since multiplication can not reduce the degree and $\Psi$ is unipotent we have
\[
\mu (\lambda_i) = \Psi \circ \pi^{h,t} (Z^M (1_0 \otimes \sml )^{-1} (D_i(Z^M ( \sml )))^{b_i})  = 1+ O(n),
\]
and the result follows since all $\mu$-invariants of length $\leq$ n vanish if and only if $\lambda_i$ is trivial in $\pi(\bsig - \sigma) / \pi(\bsig - \sigma)_n$, where $\sigma$ is a tangle represented by $\sml$.

\smallskip

Conversely, suppose that $\lambda_i$ is trivial in $\pi(\bsig - \sigma) / \pi(\bsig - \sigma)_n$.
 Then
\[\coeff{x_{\iota_1} x_{\iota_2} \cdots x_{\iota_{n-r-1}} }{\mu (\lambda_i )} = 0  \text{, \: for } 0 \leq r \leq n-1,\]
 and so
\[
\mu (\lambda_i) = \Psi \circ \pi^{h,t} (Z^M (1_0 \otimes \sml )^{-1} (D_i(Z^M ( \sml )))^{b_i})  = 1+ O(n).
\]
As $\Psi$  is unipotent it follows that
\[
\pi^{h,t} (Z^M (1_0 \otimes \sml )^{-1} (D_i(Z^M ( \sml )))^{b_i})  = 1+ O(n).
\]
Thus $D_i(\zmht{\sml}) = 1 + O(n)$ and so $\zmt{\sml} = 1 + O(n)$.

\medskip

\noindent (ii)
Suppose that the first non-vanishing $\mu$-invariant is of degree $n$. Then by the above
$\zmht{\sml} = 1 +\xi +O(n+1)$, where $\xi$ is of degree $n$.
By proposition~\ref{mu:mu}, lemma~\ref{mu:double} and the unipotency of  $\Psi$  we have
\[
\mu (\lambda_i) = \Zmht_{\leq n}(1_0 \otimes \sml )^{-1}
\cdot \Zmht_{\leq n} (\beta_i)
\cdot \Zmht_{\leq n} ( D_i(\sml ))
\cdot \Zmht_{\leq n} (\beta_i)^{-1} + O(n+1).
\]
Since $\zm{\beta_i}$ can be written as the exponential  of a sum of connected elements, this can be written as
\[
\begin{split}
&=(1 - 1_0 \otimes \xi)(1+\zeta)(1+ \pi^{h,t} D_i (\xi))(1-\zeta) + O(n+1) \\
&= 1+ \pi^{h,t} D_i (\xi)+ O(n+1)
\end{split}
\]
where $\zeta = \Zmht_n (\beta_i)$ and $\xi$ is as in the statement of the theorem.

Looking at the degree $n$ part of this formula we see that the terms of $ D_i (\xi)$ which do not have a vertex on the 0-coloured skeleton component cancel with the terms of $1_0 \otimes \xi$, and any terms of $ D_i (\xi)$ with more than one vertex on the 0-coloured skeleton component are killed off by the projection $\pi^{h,t}$.
So what remains is an element of $\mathcal{C}^t_n(X_L, 0) = \text{Lie}_n (l)$ and it is easy to see that this is exactly the element $j_i (\xi )$.

Finally, the determined by part follows since $j_i$ is injective (see \cite{HM:00}).
\end{proof}

The struts of the Kontsevich integral of a link in $S^3$ determine and are determined by its linking numbers. Theorem~\ref{mu:formulae} gives the analogous result for links in homology spheres.
\begin{Cor}
The coefficients of the struts of the \AA rhus integral of a link in a homology sphere determine and are determined by the linking numbers.
\end{Cor}

\chapter{A Diagrammatic formula for the Free Energy} \label{chapter:exp}
 
\def\bgap{\vspace{4cm}}
\def\coeff#1#2{\text{coeff}(#1,#2) }
\def\sgap{\vspace{0.5cm}}
\def\lmo#1{\langle #1 \rangle }
\def\clmo#1{ \langle #1 \rangle_c }
\def\rlmo#1#2{\langle #1 | #2 \rangle}
\def\SUM#1#2{\sum_{#1}^{#2}}
\def\graphs{\prod_{i =1}^{m} \Gamma_{i}^{p_i}}

\def\bra#1#2{\langle #1 , #2 \rangle}
\def\struts{\SUM{i=1}{n} a_i \hspace{2mm} y_i \frown y_i    }

\def\con{\SUM{r=1}{m} \Lambda_{r} }
\def\cbra#1#2{\langle #1 , #2 \rangle_c}

\def\LMO#1{ \left\langle\begin{array}{c}#1 \end{array}\right\rangle }
\def\N0{\mathbb{N}_0}
\def\bigbrac#1{\left(\begin{array}{c}#1 \end{array}\right)}
\def\CLMO#1{ \left\langle\begin{array}{c}#1 \end{array}\right\rangle_{c} }
\def\BRA#1#2{ \left\langle\begin{array}{c}#1 , #2 \end{array}\right\rangle }
\def\CBRA#1#2{ \left\langle\begin{array}{c}#1 , #2 \end{array}\right\rangle_{c} }

We state and prove a folklore result concerning the diagrammatic integration of exponentials.

\section{Introduction}

As we have already seen in chapter~\ref{chapter:mu}, the notions of diagrammatic or Feynman integration and diagrammatic differential operators play an important role in quantum topology,  for example they are used to defining  finite-type 3-manifold invariants from the Kontsevich integral (\cite{{BGRT:AI},{BGRT:AII},{BGRT:AIII}}) and certain vectorspace isomorphisms in the wheeling theorem (\cite{BLT:03}).

In physics there is a folklore principle which says that the diagrammatic integration of the exponential of something connected is itself an exponential.  We give an exact formulation of this statement.
The resulting formula is a useful combinatorial identity which can be used to simplify calculations. In particular there are applications to  the LMO invariant of 3-manifolds.

Although motivated by the theory of finite-type invariants, we find it convenient to work in a slightly more general setting in this chapter.

The version of the proof of theorem~\ref{theorem1} presented here was suggested by Daan Krammer which greatly improved an earlier proof by the author.

\begin{Remark}
We note that the author is currently  collaborating with D.M.Jackson and A. Morales to give a more classical combinatorial exposition of the diagrammatic integration discussed here and in chapter~1.
\end{Remark}

\section{Statement of Results}

Let $\mathcal{D}(Y)$ be the algebra of formal power series of
 uni-trivalent graphs  with coefficients in  $\mathbb{Q}$,  whose uni-valent vertices are coloured by some set $Y$ and trivalent vertices are oriented  and  where commutative multiplication is given by disjoint union. Note that we allow the empty graph $\emptyset$.  We also allow $Y=\emptyset$, in which case the formal power series in $\mathcal{D}(\emptyset)$ contain only  trivalent graphs.

Recall, an element of $\mathcal{D}(Y)$ is called  connected if all of its summands consist of connected graphs.

We say that $D \in \mathcal{D} (Y)$ is {\em Y-substantial} if it contains no struts (graphs which look like $\frown$). We will denote the subalgebra of $Y$-substantial elements of $\mathcal{D}(Y)$ by $\mathcal{D}_s (Y)$.

\begin{figure}
\[
\left\langle \begin{array}{c}
 \raisebox{0.0cm}{$\frac{1}{2}$}
 \hspace{1mm}
 \raisebox{0.0cm}{\epsfig{file=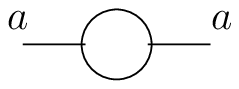, height=.3cm} }
\raisebox{0.0cm}{$+$}
\hspace{1mm}
\raisebox{0.0cm}{\epsfig{file=expfigures/aawheel2, height=.3cm}}
\raisebox{-0.2cm}{\epsfig{file=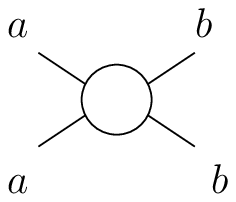, height=.7cm}}
 \end{array} \right\rangle
=
 \raisebox{0.0cm}{$\frac{1}{2}$}
\hspace{1mm}
\raisebox{-0.1cm}{\epsfig{file=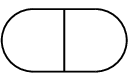, height=.3cm}}^2
 \hspace{1mm}
 \raisebox{0.0cm}{+}
\hspace{1mm}
\raisebox{-0.1cm}{\epsfig{file=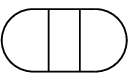, height=.3cm}}
 \hspace{1mm}
 \raisebox{0.0cm}{+ 2}
\hspace{1mm}
\raisebox{-0.1cm}{\epsfig{file=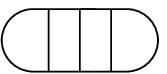, height=.3cm}}
 \hspace{1mm}
 \raisebox{0.0cm}{$+$}
\hspace{1mm}
\raisebox{-0.1cm}{\epsfig{file=expfigures/theta.eps, height=.3cm}}
 \hspace{1mm}
\raisebox{-0.1cm}{\epsfig{file=expfigures/2theta.eps, height=.3cm}}
\]
\[
\left\langle \begin{array}{c}
 \raisebox{0.0cm}{$\frac{1}{2}$}
 \hspace{1mm}
 \raisebox{0.0cm}{\epsfig{file=expfigures/aawheel2.eps, height=.3cm} }
\raisebox{0.0cm}{$+$}
\hspace{1mm}
\raisebox{0.0cm}{\epsfig{file=expfigures/aawheel2, height=.3cm}}
\raisebox{-0.2cm}{\epsfig{file=expfigures/aabbwheel4, height=.7cm}}
 \end{array} \right\rangle_c
=
\raisebox{-0.1cm}{\epsfig{file=expfigures/2theta.eps, height=.3cm}}
 \hspace{1mm}
 \raisebox{0.0cm}{+ 2}
\hspace{1mm}
\raisebox{-0.1cm}{\epsfig{file=expfigures/3theta.eps, height=.3cm}}
 \]
 \caption{Examples of definition~\ref{def}.}
\label{fig:def1}
\end{figure}

\begin{Def}\label{def}
Let $D \in \mathcal{D}_s(Y) $ be Y-substantial.  $\lmo{D}$ is defined to be the linear operation given by
\[ \lmo{D} =
\left( \begin{array}{l}
 \text{sum of all ways of identifying pairwise} \\
\text{all of the } y \text{ coloured uni-valent vertices} \\
\text{of }  D \text{ for all } y \in Y
\end{array} \right).
\]
This sum is declared to be zero if D has an odd number of $y$-coloured vertices for any $y \in Y$.
\noindent
Further, define $\clmo{D}$ to be  the connected part of $\lmo{D}$.
See figure~\ref{fig:def1} for examples of these definitions.
\end{Def}

\begin{Remark}If we work with the algebra $\mathcal{B}(Y)$ from the theory of finite-type invariants (see section~\ref{sec:algebras} or \cite{BN:95:2}), restrict ourselves to those elements with exactly $2m$ legs of each colour and project the result of $\lmo{-}$ onto the quotient of $\mathcal{B}(Y)$ by the
 the ideal generated by the relations $O_m$ and $P_{m+1}$ defined in ~\cite{LMO}, then the above definition is  {\em negative dimensional integration}, $\int^{(m)}$,  defined in~\cite{BGRT:AIII} to give  a construction of the LMO invariant -  a universal peturbative invariant of rational homology spheres.
\end{Remark}

The relationship stated in the following theorem was conjectured independently by the author and
Stavros Garoufalidis.
The proof of the theorem is given in section~\ref{proofs}.

\begin{The} \label{theorem1}
Let $C \in \mathcal{D}_s (Y) $ be Y-substantial and contain only connected graphs.  Then
\begin{equation} \label{th1}
\lmo{\exp (C)} = \exp  \left( \SUM{j=1}{\infty} \frac{1}{j!}  \left\langle C^j   \right\rangle_c  \right) .
\end{equation}
\end{The}

\sgap

Of course $\lmo{-}$ is not the only type of diagrammatic integration  in the literature and so it is natural to ask which of them satisfy equations analogous to equation~\ref{th1}.

There is a well known bilinear pairing $\bra{-}{-} : \mathcal{D} (Y) \otimes \mathcal{D} (Y) \rightarrow  \mathcal{D} (\emptyset )$, where $\bra{D_1}{D_2}$ is
 defined to be the sum of all ways of gluing all of the $y$-coloured vertices of $D_1$ to all of the $y$-coloured vertices of $D_2$ for every $y \in Y$, where this sum is zero if the number of $y$-coloured vertices of $D_1$ and $D_2$ do not match.
See also section~\ref{sec:aarhus}.
We also define $\cbra{-}{-}$ to be the connected part of $\bra{-}{-}$.

For motivation we give two important examples of uses of this definition. The first example is $\lmo{D} = \bra{ \exp ( \sum_{y \in Y} \frac{1}{2}\hspace{1mm}(y \frown y))}{D}  $ which relates our present discussion to our earlier discussion and to the LMO invariant.

Our second important example is essentially the gluing from  formal Gaussian integration (defined in section~\ref{sec:aarhus}) which is of the form
\[
\BRA{ \exp \bigbrac{ \sum_{1 \leq i \leq j \leq n} a_{i,j} \;  y_i \frown y_j } }{ D },
\]
where $D$ is Y-substantial.

The following generalizes theorem~\ref{theorem1}.
 Again we delay the proof until section~\ref{proofs}.

\begin{The} \label{theorem2}
Let $Y = \{ y_1, \cdots , y_n \} $, $C \in \mathcal{D}_s (Y)$ be $Y$-substantial and  connected and $a_{i,j} \in \mathbb{Q}$.
 Then the following identity holds
\begin{multline} \label{th2}
\left\langle \exp \left( \sum_{1 \leq i \leq j \leq n} a_{i,j} \hspace{1.5mm} y_i \frown y_j  \right)   , \exp ( C ) \right\rangle
\\
=
\exp \left(  \sum_{p=1}^{\infty} \frac{1}{p!} \left\langle \exp \left( \sum_{1 \leq i \leq j \leq n} a_{i,j} \hspace{1.5mm}  y_i \frown y_j  \right) , C^p \right\rangle_{\! \! \! c} \hspace{1mm}\right).
\end{multline}
\end{The}

\sgap
So far we have only discussed the free algebra $\mathcal{D} (Y)$.  However as we are usually interested in a quotient of $\mathcal{D} (Y)$ by some ideal we observe the following corollary, which relates the theorems presented here to the theory of finite type invariants.
\begin{Cor}
Let  $\mathcal{B} (Y)$ be the algebra of formal power series of uni-trivalent graphs with oriented trivalent vertices and univalent vertices coloured by $Y$, modulo the {\bf IHX} and {\bf AS} relations (see figure~\ref{fig:relations} for their definitions). Let $C \in \mathcal{B}(Y)$ contain only connected graphs and be $Y$-substantial.
Then equation~\ref{th2} holds in this quotient algebra.
\end{Cor}

\section{The Proofs} \label{proofs}


For motivation, consider the calculation of some value $\langle C \rangle$.  One approach is to sum over all of the ways of breaking the computation down into the construction of connected components, for example
\begin{multline*}
\langle \epsfig{file=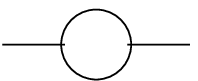, height=2mm}^2
\hspace{1mm}
 \epsfig{file=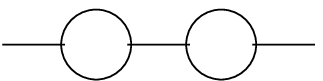, height=.2cm} \rangle
=
\langle \epsfig{file=expfigures/wheel2.eps, height=2mm}^2
\hspace{1mm}
 \epsfig{file=expfigures/glasses, height=.2cm} \rangle_c
+
\langle \epsfig{file=expfigures/wheel2.eps, height=2mm}^2 \rangle_c
\langle \epsfig{file=expfigures/glasses, height=.2cm} \rangle_c
+
2\langle \epsfig{file=expfigures/wheel2.eps, height=2mm}
\hspace{1mm}
 \epsfig{file=expfigures/glasses, height=.2cm} \rangle_c
\langle \epsfig{file=expfigures/wheel2.eps, height=2mm} \rangle_c
\\
+
\langle \epsfig{file=expfigures/wheel2.eps, height=2mm} \rangle_c
\langle \epsfig{file=expfigures/wheel2.eps, height=2mm} \rangle_c
\langle \epsfig{file=expfigures/glasses, height=.2cm} \rangle_c .
\end{multline*}
Splitting the calculation like this is the idea behind the proof of theorem~\ref{theorem1}.

More precisely the idea of the proof is to rewrite $\langle \exp (C) \rangle$ in terms of $\clmo{-}$ as indicated above, and after some rearranging of the indexing, the result drops out.
The first problem encountered is how to enumerate this sum.  We taken care of this first.

\begin{Lem}\label{sumovcon}
Let $C_i$ be a connected element of $\mathcal{D}_s (Y)$ for each $i \in I$, where  $I$ is an index. Then
\begin{equation*} 
\LMO{ \prod_{i\in I} C_i^{h(i)}}
=
\sum_{m=1}^{\infty}
\sum_{\substack{\{ g_j :I \rightarrow \N0 |   \\
\sum_j g_j =h ,\\
 g_j \neq 0, \\
1 \leq j \leq m \} }}
\frac{1}{m!} \prod_{i \in I} \frac{h(i)!}{\prod_{j=1}^{m} g_j(i)! }
 \CLMO{\prod_{i \in I} C_i^{g_j(i)}},
\end{equation*}
where $h:I \rightarrow \N0$ is a map such that only finitely many of the $h(i)$ are non-zero.
\end{Lem}

We introduce some notation for the proof.  Let $C = \prod_{i \in I} C_i^{g(i)}$ be an element of $\mathcal{D} (Y)$ such that each $C_i$ is connected.
By a {\em pattern of order $m$} of $C$ we mean a way of partitioning $C$ into
$ \prod_{j=1}^{m} \clmo{ \prod_{i \in I } C_i^{h(i,j)}} $,
where $ \sum_{ j=1 }^m h(i,j) = g(i) $ for each $i$.
Note there is no ordering of the $C_i$ or the $\clmo{-}$.
For example, with $C$ as in the above example,
$   \langle \epsfig{file=expfigures/wheel2.eps, height=2mm}
\hspace{1mm}
 \epsfig{file=expfigures/glasses, height=.2cm} \rangle_c
\langle \epsfig{file=expfigures/wheel2.eps, height=2mm} \rangle_c
 $
and
$
\langle \epsfig{file=expfigures/wheel2.eps, height=2mm} \rangle_c
\langle \epsfig{file=expfigures/wheel2.eps, height=2mm} \rangle_c
\langle \epsfig{file=expfigures/glasses, height=.2cm} \rangle_c
$
are patterns of order 2 and 3 respectively.

We now prove the lemma.
\begin{proof}
Clearly
\begin{equation} \label{patternlemmaproof}
\sum_{\substack{\{ g_j :I \rightarrow \N0 |   \\
\sum_j g_j =h ,\\
 g_j \neq 0, \\
1 \leq j \leq m \} }}
\CLMO{\prod_{i \in I} C_i^{g_j(i)}}
\end{equation}
enumerates all patterns of order $m$. It remains to add  coefficients to count the correct number of occurrences of each pattern to equation~\ref{patternlemmaproof} which, after summing over all values of  $m$, will give the equation in the lemma.

As each pattern in equation~\ref{patternlemmaproof} occurs with a multiplicity of $m!$ (as the sum orders the product of $\clmo{-}$ in the pattern) we need to divide each summand by $m!$ to give exactly one occurrence of each pattern. The lemma then follows upon noting that there are
\[
\prod_{i \in I}
\binom{h(i)}{g_1(i), \cdots , g_m(i)}
=
\prod_{i \in I} \frac{h(i)!}{\prod_{j=1}^{m} g_j(i)! }
\]
ways of making the pattern
$ \CLMO{\prod_{i \in I} C_i^{g_j(i)}} $
from
$ \prod_{i\in I} C_i^{h(i)} $.
\end{proof}

\medskip

\begin{proof}[Proof of theorem~\ref{theorem1}]
We begin by writing $C$ as $\sum_{i \in I} C_i$, where each $C_i$ is a connected element of $\mathcal{D}_s (Y)$ and $I$ is an indexing set. This gives
\[
\lmo{\exp (C)}
= \LMO{ \prod_{i \in I} \exp (C_i)}
= \LMO{\prod_{i \in I} \sum_{n=1}^{\infty} \frac{C_i^n}{n!}  }.
\]
This can be written as a sum over functions in the following way:
\[
\LMO{\sum_{h:I \rightarrow \N0} \prod_{i \in I} \frac{C_i^{h(i)}}{h(i)!} }
=
\sum_{h:I \rightarrow \N0} \LMO{\prod_{i \in I} \frac{C_i^{h(i)}}{h(i)!}}.
\]

Now, using lemma~\ref{sumovcon}, we can  split the calculation of $\lmo{-}$ into a sum of products of $\clmo{-}$, giving
\[
1+ \sum_{h:I \rightarrow \N0}
 \sum_{m=1}^{\infty} \frac{1}{m!}
\sum_{\substack{\{ g_j :I \rightarrow \N0 |   \\
\sum_j g_j =h,\\
 g_j \neq 0, \\
1 \leq j \leq m \} }}
\bigbrac{ \prod_{i \in I}  \frac{1}{h(i)!}   \frac{h(i)!}{ \prod_{j=1}^{m} g_j (i)!  }  }
\prod_{j=1}^{m}
\CLMO{\prod_{i \in I} C_i^{g_j (i)}}.
\]
We may incorporate the sum over the functions $h$ into the sum over the $g_j$ to write this as
\[
1+
 \sum_{m=1}^{\infty} \frac{1}{m!}
\sum_{\substack{\{ g_j :I \rightarrow \N0 |   \\
 g_j \neq 0, \\
1 \leq j \leq m \} }}
\prod_{j=1}^{m}
\CLMO{\prod_{i \in I} \frac{C_i^{g_j (i)}}{g_j (i)!}},
\]
which with a little thought,  can be seen to be equal to
\[
1+
 \sum_{m=1}^{\infty} \frac{1}{m!}
\prod_{j=1}^{m}
\bigbrac{
\sum_{\substack{ g_j :I \rightarrow \N0 ,   \\ g_j \neq 0  }}
\CLMO{\prod_{i \in I} \frac{C_i^{g_j (i)}}{g_j (i)!}}},
\]
where  the sum is now over $g_j$ for a fixed $j$.
Clearly this may be written as
\[
 \sum_{m=0}^{\infty} \frac{1}{m!}
\bigbrac{
\sum_{\substack{ g :I \rightarrow \N0 ,   \\ g \neq 0  }}
\CLMO{\prod_{i \in I} \frac{C_i^{g (i)}}{g (i)!}}}^m
\]
\[
=\exp
\bigbrac{
\sum_{\substack{ g :I \rightarrow \N0 ,   \\ g \neq 0  }}
\CLMO{\prod_{i \in I} \frac{C_i^{g (i)}}{g (i)!}}}
\]
It remains to show that
\[
\sum_{\substack{ g :I \rightarrow \N0 ,    \\ g \neq 0  }}
\CLMO{\prod_{i \in I} \frac{C_i^{g (i)}}{g (i)!}}
=
\sum_{j=1}^{\infty} \frac{1}{j!} \CLMO{C^j}
\]
but this follows as one may write the right hand side as
\[
\sum_{j=1}^{\infty} \frac{1}{j!} \sum_{\substack{ g :I \rightarrow \N0 ,   \\ \sum_{i \in I} g(i) =j   }}
\frac{j!}{\prod_{i \in I} g(i)!}
\CLMO{\prod_{i \in I} C_i^{g(i)} },
\]
which is obviously equal to the left hand side.
\end{proof}

\sgap
Since the proof of theorem~\ref{theorem2} is similar to that of theorem~\ref{theorem1} we only sketch it.
The only real difference is that in a calculation of $\bra{S}{C}$ in terms of $\cbra{-}{-}$, we look at all the splittings of both  $S$ and $C$, for example
\begin{multline*}
\langle  \frown^{3} ,
\epsfig{file=expfigures/wheel2.eps, height=2mm}^2
\hspace{1mm}
 \epsfig{file=expfigures/glasses, height=.2cm} \rangle
=
\langle
\frown^{3} ,
\epsfig{file=expfigures/wheel2.eps, height=2mm}^2
\hspace{1mm}
 \epsfig{file=expfigures/glasses, height=.2cm} \rangle_c
+
\langle \frown^{2} , \epsfig{file=expfigures/wheel2.eps, height=2mm}^2 \rangle_c
\langle \frown , \epsfig{file=expfigures/glasses, height=.2cm} \rangle_c
\\
+
\langle \frown , \epsfig{file=expfigures/wheel2.eps, height=2mm}^2 \rangle_c
\langle \frown^{2} , \epsfig{file=expfigures/glasses, height=.2cm} \rangle_c
+
\cdots
+
\langle \frown , \epsfig{file=expfigures/wheel2.eps, height=2mm} \rangle_c
\langle \frown , \epsfig{file=expfigures/wheel2.eps, height=2mm} \rangle_c
\langle \frown , \epsfig{file=expfigures/glasses, height=.2cm} \rangle_c .
\end{multline*}
In fact the proof presented below is essentially a double application (to the struts and to $C$) of the proof of theorem~\ref{theorem1}.

\begin{proof}[Sketch of the proof of theorem~\ref{theorem2}]
As before we write $C$ as $\sum_{i \in I} C_i$ where each $C_i$ is a connected element of $\mathcal{D}_s (Y)$ and $I$ is an indexing set. Also write the strut part as $\sum_{k \in K} S_k$, for some index $K$.  Then

\[
\BRA{\exp \left( \sum_{k \in K} S_k \right)}{\exp (C)}
= \BRA{\prod_{k \in K} \sum_{m=1}^{\infty} \frac{S_k^m}{m!}}{\prod_{i \in I} \sum_{n=1}^{\infty} \frac{C_i^n}{n!}  }.
\]
As before this can be written as a sum over functions:
\[
\sum_{h:I \rightarrow \N0}
\sum_{\alpha :K \rightarrow \N0}
\BRA{\prod_{k \in K} \frac{S_k^{\alpha(k)}}{\alpha(k)!} }{\prod_{i \in I} \frac{C_i^{h(i)}}{h(i)!}}.
\]
By using an argument similar to that of the proof of lemma~\ref{sumovcon}, we write this as
\begin{multline*}
1+ \sum_{h:I \rightarrow \N0}
 \sum_{m=1}^{\infty} \frac{1}{m!}
\sum_{\substack{\{ g_j :I \rightarrow \N0 |   \\
\sum_j g_j =h,\\
 g_j \neq 0, \\
1 \leq j \leq m \} }}
\bigbrac{ \prod_{i \in I}  \frac{1}{h(i)!}   \frac{h(i)!}{ \prod_{j=1}^{m} g_j (i)!  }  }
\\
\sum_{\alpha :I \rightarrow \N0}
 \sum_{n=1}^{\infty} \frac{1}{n!}
\sum_{\substack{\{ f_l :K \rightarrow \N0 |   \\
\sum_l f_l =\alpha,\\
 f_l \neq 0, \\
1 \leq l \leq n \} }}
\bigbrac{ \prod_{k \in K}  \frac{1}{\alpha(k)!}   \frac{\alpha(k)!}{ \prod_{l=1}^{m} g_l (k)!  }  }
\\
\prod_{j=1}^{m}
\prod_{l=1}^{n}
\CBRA{ \prod_{k \in K} S_k^{f_l (k)}  }{\prod_{i \in I} C_i^{g_j (i)}}
\end{multline*}
\begin{multline*}
=1+ \sum_{h:I \rightarrow \N0}
 \sum_{m=1}^{\infty} \frac{1}{m!}
\sum_{\substack{\{ g_j :I \rightarrow \N0 |   \\
\sum_j g_j =h,\\
 g_j \neq 0, \\
1 \leq j \leq m \} }}
\bigbrac{ \prod_{i \in I}  \frac{1}{h(i)!}   \frac{h(i)!}{ \prod_{j=1}^{m} g_j (i)!  }  }
\\
\prod_{j=1}^{m}
\CBRA{ \exp ( \sum_{k \in K} S_k )  }{\prod_{i \in I} C_i^{g_j (i)}}
\end{multline*}
\[
=
\exp \bigbrac{  \sum_{p=1}^{\infty} \frac{1}{p!} \CBRA{\exp \left( \sum_{k \in K} S_k \right) }{ C^p }    },
\]
using similar indexing arguments to those above.
\end{proof}

\begin{Remark}
An earlier proof of theorem~\ref{theorem2} gave a partial converse to the theorem, however as it has been superseded by results to appear in \cite{JMcomm}, we do not include it.
\end{Remark}

\chapter{On the Word and Conjugacy Problems for Link Groups} \label{chapter:word}

\def\bgap{\vspace{4cm}}
\def\sgap{\medskip}

\def\canc4{$C^{\prime \prime} (4)$}
\def\cant4{$T (4)$}
\def\c4t4{$C^{\prime \prime} (4) - T (4)$ }
\def\eqg{=_G}
\def\eqpi{=_{\pi}}
\def\neqg{\neq_G}
\def\neqpi{\neq_{\pi}}
\def\vec#1{\mathbf{#1}}
\def\adp{augmented Dehn presentation }
\def\pres#1#2{\langle #1 | #2 \rangle}

In this chapter we study the relationship between known  solutions from small cancellation theory and normal surface theory for the word  and conjugacy problems of the groups of (prime) alternating links.
We show that two of the algorithms in the literature for solving the word problem, each using one of the two approaches, are the same.  Then, by considering small cancellation methods, we give a normal surface solution to the conjugacy problem  of these link groups and characterize the conjugacy classes.
Finally, by using the small cancellation properties of link groups we provide a new proof that alternating links are non-trivial.

\section{Introduction}
The word problem for the fundamental group of a link complement was first solved by Waldhausen in \cite{Wa1968}, who found an algorithm for deciding whether a loop in a sufficiently large irreducible 3-manifold is contractible.
Waldhausen's algorithm is difficult to apply. It relies on Haken's theory of normal surfaces to find a particular set of surfaces in the 3-manifold.
Dugopolski, in \cite{Du1982}, showed that for the complement of an alternating link, such a set of surfaces is readily available and he simplified Waldhausen's solution of the word problem for such links.  These algorithms are geometric.

On the other hand, using combinatorial group theory, Weinbaum, in \cite{We1971},  proved that the  groups of prime alternating knots satisfy the  \c4t4 small cancellation conditions, and applied results of Lyndon, to solve the word problem, and Schupp, to solve the conjugacy problem.
 This was extended to all alternating knots by Appel and Schupp in \cite{AS1972} (see also \cite{LSbook}) by showing these groups satisfy the \c4t4 for minimal sequences conditions.
In fact these methods apply to a larger set of (not necessarily alternating) links.

At around the same time  Appel, \cite{Ap1974},  used the Wirtinger presentation and small cancellation techniques to solve the conjugacy problem for all alternating and some non-alternating knots.

 An improvement for solving the word problem for \c4t4 groups and hence the groups of prime alternating links,  to a quadratic time algorithm, appears implicitly in later work by Appel and Schupp (\cite{AS1983}), and explicitly in several other places (\cite{{GS1990},{Jothesis},{Jo2000},{Ka1997}}). We  consider  this algorithm.

Other approaches have also been successful.  Epstein and Thurston in \cite{Epbook} showed that all link groups are automatic (Gersten and Short also showed that the groups of alternating knots are automatic using different methods in \cite{GS1991}) and hence have a solvable word problem.
In fact the conjugacy problem for link groups is solvable in full generality since the link complement is a CAT(0) space.

Returning to small cancellation methods, Johnsgard gave a polynomial time algorithm for solving the conjugacy problem for prime alternating link groups in \cite{Jo1997}.

\smallskip

In this chapter we show that the Dugopolski's solution of the word problem, which uses the theory of normal surfaces, and the algorithm which comes from the characterization of geodesics in \c4t4 small cancellation groups, are  the same for prime alternating links.

More specifically, we see that, expressed in terms of \c4t4 groups, the two algorithms differ only in the way that they search for subwords to freely reduce and carry out chain collapses (a particular type of substitution), also Dugopolski's algorithm will cyclically permute the word, while the other algorithm does not. These differences are inconsequential.  In terms of the language of Dugopolski, the algorithm from group theory carries out type 1 reductions and type 2 deformations with respect to both the white and black checker-board hierarchy, while keeping one point of the loop fixed, while Dugopolski's carries out the moves with respect to one hierarchy only and does not fix a  point.

We go on to use the correspondences developed in  proving this  to give a normal surfaces algorithm for solving the conjugacy problem for prime alternating links.
This is important since normal surface theory has not been particularly successful in solving the conjugacy problem, although Evans (\cite{Ev1978}) used the theory to solve the conjugacy problem for loops in the boundary of a compact sufficiently large 3-manifold.
 This also provides a geometric characterization of the conjugacy classes of a given loop in the link complement.
We also give a characterization of the classes of contractible loops using these moves.

A secondary purpose of this chapter is to provide a readable account of the methods and algorithms used herein.
We  discuss in detail the processes needed to apply Dugopolski's algorithm, which was not discussed in his paper and we also give a short proof of the planarity of Johnsgard's conjugacy algorithm in the difficult case (which is needed to prove that it is polynomial time).

Finally we use the small cancellation properties of link groups to give a new proof of the non-triviality of alternating links.

\begin{Remark}
Although we will mostly talk about prime links, one should remember that the word and conjugacy problem for the free product of two groups reduces to that of its factors, so we obtain solutions for the split unions of the links.  Also note that the free product of \c4t4 groups is again \c4t4 and so the small cancellation results also hold for split unions.
\end{Remark}

\section{Definitions}

\begin{figure}
\centering
\epsfig{file=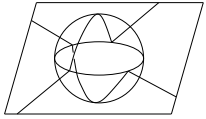, height=1.5cm}
\caption{A crossing ball.}
\label{crossingballs}
\end{figure}

Let $L \subset \mathbb{R}^3 \cup \infty = S^3$ be a link.
 Without loss of generality we may assume that $L$ lies on $ \mathbb{R}^2 \cup \infty$ except in a neighborhood of a crossing where the arcs lie on the boundary of a 3-ball, which we call a {\em crossing ball},  forming semicircular over and under crossing arcs which intersect the north and south poles.  This is  indicated in figure~\ref{crossingballs}.  For convenience we will always assume that our links are of this form.

There is a {\em canonical projection} associated with every such link, given by vertical projection onto $\mathbb{R}^2 \cup \infty$ inside the crossing balls.  Similarly we can associate a link to every projection by adding crossing balls at each crossing.

Since we can move canonically between a link and a projection, we will abuse notation and say that a link has some property when its canonical projection does.

\begin{figure} \label{2-cells}
\centering
\epsfig{file=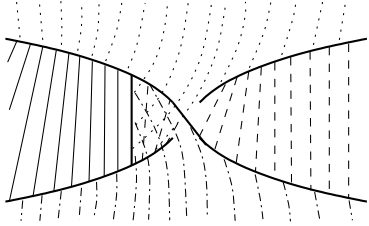, height=3cm}
\caption{Regions of a link.}
\end{figure}

We can use the regions of the canonical projection to induce a set of surfaces in the link complement. These are the 2-cells which coincide with the regions of the projection outside of the crossing balls and a strip with a $\pi /2$ twist  inside the crossing balls whose boundary is identified with the arc given by the intersection of the crossing ball and the region of the projection the region, the two arcs from the equator to the poles and the north-south axis of the crossing balls, shown locally in figure~\ref{2-cells}. We call these 2-cells the {\em regions} of the link.

\begin{figure}\label{relators}
\centering
\epsfig{file=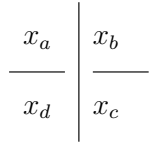, height=2cm}
\caption{Relators for the presentation.}
\end{figure}

 The {\em Dehn presentation} of $\pi_1 (S^3 -L)$ is defined as follows: Take a regular projection of $L$ onto $ \mathbb{R}^2 \cup \infty $ and label the regions  $x_0, x_1, \ldots , x_n$ (note that this induces a labeling of the regions of the link).
These will be the names of our generators.
 By convention we label the region containing infinity $x_0$, and call this the {\em outer region}.
To each crossing we assign a relator $x_a x_b^{-1} x_c x_d^{-1}$ according to figure~\ref{relators} and add one extra relator $x_0$.  The presentation obtained after we kill off the generator $x_0$ using Titze transformations is called the {\em Dehn presentation} of $L$.

A geometric interpretation of this presentation of the fundamental group  follows by
choosing a base point above $\mathbb{R}^2 \cup \infty$ and to each region $x_i$ of the link, assign a loop which descends from the base point through the region $x_i$ and returns through $x_0$.

  We say that a  set of relators is  {\em symmetrized} if it is cyclically reduced and closed with respect to inverses and cyclic permutations.  A symmetrized presentation can be obtained from any presentation by adjoining the cyclic permutations of the set of relators and their inverses.  Clearly this does not change the group.

We denote the unknot by $\mathcal{O}$.
The {\em augmented link} of $L$ is the link corresponding to the projection of $L \cup \mathcal{O}$  such that the projection of $L$ shares no edges with the outer region (so the projection of $L$ is ``inside'' the projection of $\mathcal{O}$).

By convention we label the region which is  bounded by the unknot component and contains the projection of $L$ (so not the outer region) $x_0$.

Now define the {\em augmented Dehn presentation} of $L$ to be the symmetrized Dehn presentation of the augmented link.

The augmented Dehn presentation is obtained from the construction of the Dehn presentation by failing to add the relator $x_0$ and then symmetrizing.  It is  the  free product of the Dehn presentation of $L$ and the infinite cyclic group.
Therefore, solving the conjugacy (and word) problem for the augmented link group solves it for the link group.

The inclusion  $i: S^3 -L \hookrightarrow  S^3 - (L \cup \mathcal{O})$, induces the homomorphism $i_*: x_i \mapsto x_i x_0^{-1} $ from the Dehn presentation to the augmented Dehn presentation.  We will call a word   which lies in the image of the Dehn presentation  under this induced homomorphism  an {\em included word} in the augmented Dehn presentation.

If a presentation has all relators of equal length and we can assign a parity to each generator and its inverse $x_i^{\pm 1}$ such that the letters of each relator alternate in parity, we say that the presentation has {\em parity}.

Recall that the {\em checker-board colouring} of a link projection is the assignment of the  colour black or white to each of the regions in such a way that at each crossing, adjacent regions have a different colour.
By convention we  assume that the outer region is coloured white.
It is not hard to see that the checker-board colouring induces a parity on the presentation.

\section{Small Cancellation Theory and the Word Problem}

We begin by reviewing some basic constructions from group theory.
One can associate a {\em standard 2-complex} $K$ to a group presentation $G= \langle X | R \rangle $ in the usual way: $K$ consists of one 0-cell, one labelled 1-cell for each generator and one 2-cell for each relator, where the 2-cell $D_r$ representing the relator $r \in R$ is attached to the 1-skeleton, $K^{(1)}$, by a continuous map which identifies $\partial D_r$ with a loop representing $r$ in the 1-skeleton.

A  word $w \in F(X)$ represents the identity in  $G$ if and only if there is a connected simply connected planar 2-complex $D$ and a map
$\phi : (D,\partial D) \rightarrow (K, K^{(1)})$ such that the 0-cells are mapped to  0-cells, open $i$-cells are mapped to open $i$-cells, for $i=1,2$ and $\partial D$ is mapped to the loop representing $w$.
We call such a 2-complex, labelled in the natural way, a {\em singular disc diagram} (or {\em Dehn diagram} or {\em Van Kampen diagram}).  If $D$ contains no cut vertices (those whose removal disconnects $D$) then we call $D$ a {\em disc diagram}.

We say a singular disc diagram is {\em reduced} if there are no 2-cells $R_1$ and $R_2$ with a common edge $e$ such that reading the labels on their boundaries from edge $e$ clockwise on $R_1$ and anticlockwise on $R_2$ give the same word.  It is easy to see how to remove two such 2-cells without changing the boundary word.  We assume that all singular disc diagrams are reduced.

A {\em piecewise Euclidean (PE) complex} is a combinatorial 2-complex where each 2-cell is equipped with the metric of a convex polygon in the Euclidean plane in such a way that all the metrics agree on edges common to the boundaries of more than one 2-cell.

Unless otherwise stated, in this chapter we give all of the 2-complexes a PE structure  by regarding the 2-cells as regular polygons of side 1 where the number of sides of the polygon is determined by the length of the word labelling the boundary.
Note that since we never consider presentations with relators of length 2 we can do this.

\medskip

Let $G= \langle X | R \rangle $ be a group presentation.  We call a non-empty word $r$  a {\em piece with respect to $R$} if there exist distinct words $s,t \in  R$ such that $ s = r u$ and $t = rv$ and $R$ is symmetrized. Furthermore, we say that  a symmetrized presentation is \c4t4  if it satisfies the following two small cancellation conditions:

\begin{Conc4}
All relators have length four and no defining relator is a product of fewer than four pieces.
\end{Conc4}

\begin{Cont4}
Let $r_1 , r_2$ and $r_3$ be any three defining relators such that no two of the words are inverses to each other, then one of $r_1 r_2$, $r_2 r_3$ or $r_3 r_1$ is freely reduced without cancellation.
\end{Cont4}

In this chapter we are mostly concerned with \c4t4 small cancellation groups and so we usually consider {\em square complexes}, that is PE 2-complexes where the 2-cells are regarded as solid Euclidean squares.

We call a (PE) disc diagram whose boundary is labelled by a relator of length 4 a {\em relator square}. We observe that rotating a relator square by $\pi /2$ gives a cyclic permutation of the relator and flipping the square corresponds to taking the inverse of the relator.  Also note that if we are tiling with relator squares and we have a right angle with labelled edge traversals then, by \canc4 , at most one relator square exists which can fill this angle.

We call a two letter subword of a relator a {\em pair}.  The \canc4 condition says that a pair determines a relator square uniquely up to cyclic permutation and \cant4 says that if $ab$ and $b^{-1}c$ are pairs then $ac$ is not.
It follows that in a \c4t4 group a given pair uniquely determines a second pair from the corresponding relator and these pairs are equal in the group.  We call the process of replacing one pair with the other pair it determines {\em exchanging a pair}.

\begin{figure}
\centering
\epsfig{file=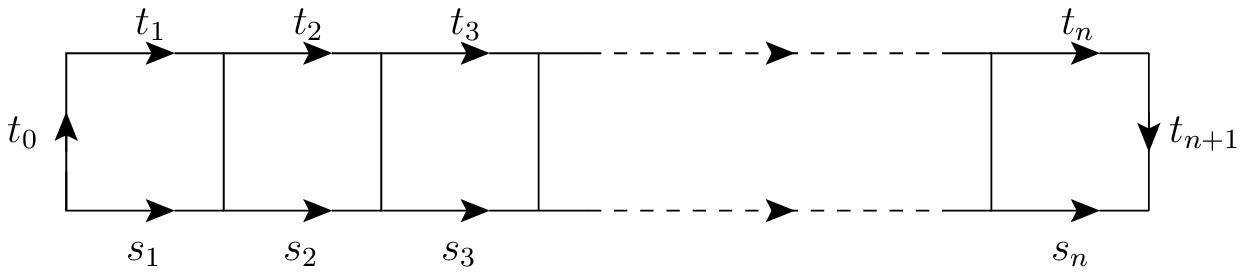, width=12cm}
\caption{A chain.}
\label{chain}
\end{figure}

A {\em chain} is a reduced $n \times 1$ disk having the form shown in figure~\ref{chain}, where $n \geq 1$. We call the word $t_0 t_1 t_2 \cdots t_{n+1}$ a {\em chain word} and the word $s_1 s_2 \cdots s_n$ the {\em inner link path}.
We collectively refer to the inner link path and the edge labelled $t_1 t_2 \cdots t_n$ as the {\em sides of the chain}.
If an inner link path of a chain is of the form $a  t_1 t_2 \cdots t_{2n}a^{-1}$ we say that it is a {\em conjugacy chain}.
 We call the process of replacing a chain word with its inner link path a {\em chain collapse}.

let $D$ be a PE disc diagram and $v$ be a 0-cell on $\partial D$. Define the {\em turning angle} to be $\tau (V) = \pi - \sigma_v $, where $\sigma_v$ is the sum of the corner angles of all of the 2-cells incident to $v$.

The following is a consequence of Lyndon's curvature theorem.

\begin{The}[\cite{GS1990}] \label{boundarychains}
Let  $G= \langle X | R \rangle $ be a \c4t4 group and $D$ be a disc diagram. Then on $\partial D$ there are at least 4 more vertices with positive turning angle than there are vertices with negative turning angle. Therefore there are at least four chains on the boundary of $D$ with disjoint interiors.
\end{The}

This immediately gives an algorithm for solving the  word problem for a word $w$ in a \c4t4 group:

\begin{Alg}\label{smallalg}
Freely reduce $w$.  If  $w$ is empty then $w \eqg 1$.
If $w$ is non-empty then scan the word from the beginning for a pair. If no pairs exist then $w \neqg 1$, otherwise scan for the next pair.  If this pair forms a chain word then perform the chain collapse and repeat all of the above.  If this pair does not form a chain word then search for the next pair and check whether this pair and the preceding pair forms a chain word.  If so perform the chain collapse and go to the beginning of the algorithm. Continue like this until all pairs have been checked.  If $w$ is still non-empty then $w \neqg 1$.
\end{Alg}

If we consider \c4t4 presentations with parity then there two types of chains: those with a white inner link path, which we call {\em white chains} and those with black inner link path which we call {\em black chains}.
The algorithm above performs chain collapses on both colours of chain. For comparison with the normal surfaces algorithm later, it is useful to know that the above algorithm works if we only carry out chain collapses on chains of a single colour. The following proposition tells us we can do this.

\begin{Prop}
Let  $G= \langle X | R \rangle $ be a \c4t4 presentation with parity and $D$ be a disc diagram. Then on $\partial D$ there exist both white chains and black chains.
\end{Prop}

\begin{proof}
To each 0-cell $v$ on the boundary $\partial D$  assign the number $3-d(v)$,  where $d(v)$ is the number of 1-cells incident to $v$. We call this number the {\em weight}.
Note that  positive (resp. zero, resp. negative) turning angles correspond to 0-cells  of weight 1 (resp. 0, resp. $<1$).
By reading off the non-zero weights, with respect to some choice of starting point, we get a sequence of numbers, which we call the {\em sequence of weights}.  We prove the result by considering the structure of this sequence and its relationship to the chains on the boundary.

First observe that the occurrence of a subsequence $1,1,1$ in the sequence of weights corresponds to a sequence of three positive turning angles. The middle turning angle must then be the corner of two chains of different colours.

Now suppose that such a sequence does not occur. Then we may assume that the sequence of weights is of the form
$ 1,1,\sigma_1 ,1,1,\sigma_2 , \cdots ,   1,1,\sigma_n $
where each $\sigma_i$ is a subsequence of weights which either consist of one negative weight or the first and last terms are negative and any positive term is bounded on each side by a negative term.
Note that these conditions imply that the sum of the weights in each $\sigma_i$ is at most $-1$.

Now if there exists a subsequence $\sigma_i$ such that the sum of its weights is odd then it is not hard to see that the two subsequences $1,1$ bounding $\sigma_i$ correspond to chains of different colours and we are done.

It remains to show that it is impossible for the sum of the weights of every  $\sigma_i$ to be even.  Suppose this was the case. Then a version of Lyndon's curvature formula (see \cite{LSbook}) gives
\[
4 \leq \sum_{v \in \partial D} (3-d(v))
= \sum_{v \in \partial D} w(v)
= 2n + \sum_{i=1}^n w(\sigma_i),
\]
where $w(v)$ is the weight of vertex $v$ and $w(\sigma_i)$ is the sum of the terms of $\sigma_i$.  But if each $w(\sigma_i) \leq -2$, this is impossible, giving the required contradiction.
\end{proof}

\medskip

Given an arbitrary finite group presentation, the set of lengths of all words representing an element of the group has a minimum. Any word which attains this minimum is called a {\em geodesic}. Geodesics in a \c4t4 presentation are characterized by the absence of chain words:
\begin{GCT}[\cite{{AS1983},{GS1990},{Jothesis},{Jo2000},{Ka1997}}]
A word in a \c4t4 presentation is geodesic if and only if it is freely reduced and contains no chain subwords.
\end{GCT}
Hence algorithm \ref{smallalg} finds a geodesic representative of a given word.  We will make use of this later.

\sgap

We say that a link projection is {\em reduced} (or {\em untwisted}) if at every crossing in the  projection, four distinct regions meet.  If in addition two distinct regions have at most one edge in common we say that the projection is {\em elementary}.  We note that every prime alternating knot has an alternating elementary projection (\cite{{LSbook},{Sc1949}}) and if an alternating projection is elementary then its corresponding link is prime (\cite{Me1984}).

\begin{The}[\cite{We1971}] \label{dehncan}
If a link has an alternating elementary projection (and is therefore  prime and alternating), then the
 augmented Dehn Presentation satisfies the \c4t4 small cancellation conditions and thus the word and conjugacy problems for its fundamental group are solvable.
\end{The}

\section{A solution to the Word Problem by Normal Surfaces}

We will now outline Dugopolski's  algorithm for deciding whether a loop in the complement of an alternating link is contractible.  Clearly  this solves the word problem for the fundamental groups of alternating links. The reader is referred to \cite{Du1982} for
the justification of the algorithm.

\sgap

\begin{figure}
\centering
\subfigure[]{\raisebox{0mm}{\epsfig{file=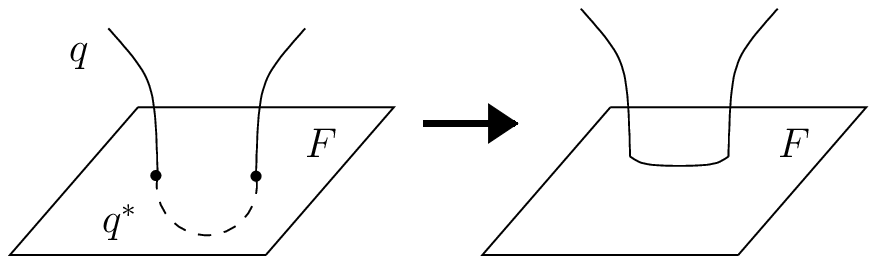, height=1.5cm}}}
\hspace{1cm}
\subfigure[]{\epsfig{file=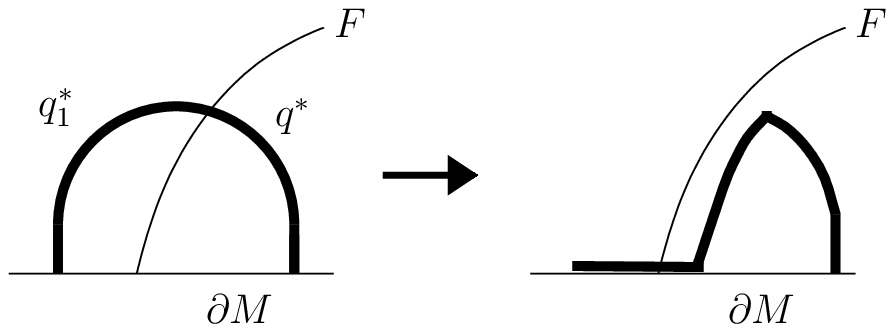, height=2cm}}
\caption{Type 1 and 2 moves.}
\label{type1and2}
\end{figure}

Let $M$ be a 3-manifold and $F$ be a surface in $M$. Suppose that  $q$ is a loop in $M$ and $q^*$ is an arc of $q$ such that $q^* \cap F = \partial q^*$. If there exists a deformation of $q^*$ into $F$ which keeps the endpoints fixed, then we may deform $q$ so that it intersects $F$ at two fewer points (see figure~\ref{type1and2}(a)). We call such a deformation {\em type 1 reduction} of $q$ with respect to $F$.

We call the inverse move (which increases the number of intersections with respect to $F$ by two) a {\em type 1 augmentation} of $q$ with respect to $F$.

Now suppose that $q^*$ is an arc of $q$ such $q^* \cap \partial M = \partial q^*$, $F$ is a surface in $M$ such that $\partial q^* \cap F = \emptyset$ and $q^*_1$ is a sub-arc of $q^*$  such that $q^*_1 \cap F$ consists of one point of  $\partial q^*_1$ and $ q^*_1 \cap \partial q^* $  is the other point of $\partial q^*_1$. If $q^*_1$ deforms to an arc lying in $F \cup \partial M$ intersecting $\partial F$ once, then $q^*$ deforms to an arc lying partly in $\partial M$, where the part not in $\partial M$ intersects $F$ one less time than $q$ does (see figure~\ref{type1and2}(b)).  We call this a {\em type 2 deformation of $q^*$}.

\sgap

\begin{figure}
\centering
\epsfig{file=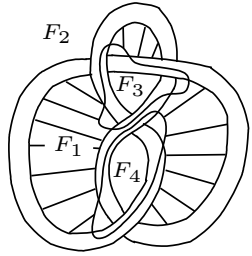, height=3.5cm}
\caption{A hierarchy for the figure 8 knot.}
\label{heirachy}
\end{figure}

We will associate a set of surfaces $F_1, \ldots F_n$ to the complement of an alternating link $L$.  Let $D \subset \mathbb{R}^2 \cup \infty \subset S^3$ be an alternating projection of $L$.
The checker-board colouring induces a colouring on the regions of the link (recall our convention of colouring the outer region white).
Form a black (resp. white) surface in the usual way by identifying all the black (resp. white) regions along their intersection at the polar axes.
This gives surfaces whose boundary is $L$ (in the case of a knot, these are spanning surfaces). Choose one of the surfaces, black say, call it $H$.  Let $F_1 = H - N(L)$, where $N(L)$ is a small regular neighborhood of $L$. Let $N(F_1)$ be a small regular neighborhood of $F_1$ and define the surfaces $F_2 , \ldots , F_n$ to be what remains of the white surfaces in $S^3 - (N(L) \cup N(F_1))$ (see figure~\ref{heirachy}).  We call the set of surfaces  $F_1 , \ldots , F_n$ the {\em black checker-board hierarchy} of $D$.
If we construct $F_1$ from the white surface, we  call the resulting set of surfaces the {\em white checker-board hierarchy}.
 Unless otherwise stated, we will use  the black checker-board hierarchy and will refer to is simply as the checker-board hierarchy.

\begin{Remark}
The hierarchy we use has one more surface (from the outer region) than the one used by Dugopolski (although the extra surface is used  implicitly at some points in his  paper).  This extra surface does not effect the validity of the algorithm.
\end{Remark}

Notice that we can split $S^3- N(L)$ into 3-balls by removing small regular neighborhoods of all the $F_i$.  In particular if we have some arc $q$ in $S^3-L$ then the only obstructions to being able to perform a type 1 reduction or a type 2 deformation with respect to some $F_i$ comes from the intersections of $q$ with the boundary of the 3-balls.  This is discussed  in more detail in section~\ref{app:petronio}, but for the time being we will just note that this gives a process to check whether we can perform a type 1 reduction or a type 2 deformation with respect to any of the $F_i$.  Retaining the notation of \cite{Du1982}, we give this the mysterious title {\em process X} with respect to $F_i$.

\sgap

Dugopolski's paper solves the word problem for the fundamental groups of certain 3-manifolds in which there exist a finite set of surfaces which satisfy certain technical conditions. He provides an algorithm which decides whether a type 1 reduction with respect to one of these surfaces is possible for a given arc.  He then goes on to prove that if a loop contracts then   the number of times that it intersects with the set of surfaces can be reduced to zero by type 1 reductions.  He  shows that the checker-board hierarchy is such a set of surfaces in the complement of alternating links and provides and algorithm for deciding whether a loop in the complement is contractible or not. The algorithms follow.

\sgap

In what follows let $L$ be an alternating link and $F_1, \ldots , F_n$ be the corresponding checker-board hierarchy. Also let  $l$ be a loop in $S^3 - N(L)$ and $l_j$ be an arc of $l$ such that $l_j \cap N(F_1) = \partial l_j$.

The following algorithm  determines whether it is possible to perform a type 1 reduction on $l_j$ with respect to $F_1$.

\begin{Alg}[\cite{Du1982}] \label{f1alg}
If $l_j \cap \cup_{i=2}^n F_i = \emptyset$ then use process X to check if a type 1 reduction with respect to $F_1$ is possible.  If $l_j \cap \cup_{i=2}^n F_i \neq \emptyset$, then for every subarc of $l_j$ with endpoints on the same disc, check for and perform type 1 reductions with respect to that surface.  Repeat this until  no more such reductions can be made.
If now $l_j \cap \cup_{i=2}^n F_i = \emptyset$ then use process X to check if a type 1 reduction with respect to $F_1$ is possible. If  $l_j \cap \cup_{i=2}^n F_1 \neq \emptyset$, then split $l_j$ into arcs $l_{j1}, \ldots , l_{js}$ by $ \cup_{i=2}^n F_i$. Use process X to check $l_{j1}$ for a type 2 deformation. If one is possible deform $l_j$ so that it lies partly in $\partial N(F_1)$ and the part which is not in    $\partial N(F_1)$ is an arc with endpoints in  $\partial N(F_1)$, intersecting  $\cup_{i=2}^n F_i$ one fewer times. Now if  $l_j \cap \cup_{i=2}^n F_i = \emptyset$ then use process X to check if a type 1 reduction with respect to $F_1$ is possible and if
$l_j \cap \cup_{i=2}^n F_i \neq \emptyset$, then repeat this step for the new arc $l_{j2}  l_{j3} \cdots l_{js}$ until no further type 2 reductions may be performed. If now $l_j \cap \cup_{i=2}^n F_i = \emptyset$ then use process X to check if a type 1 reduction with respect to $F_1$ is possible.  If
$l_j \cap \cup_{i=2}^n F_i \neq \emptyset$ then a type 1 reduction with respect to $F_1$ is not possible.
\end{Alg}

\sgap

This algorithm decides whether a type 1 reduction is possible, however to implement the algorithm we also need to know what the loop will look like after we have carried out the type 1 reductions and the type 2 deformations.
It is clear how to draw an arc after a type 1 reduction.
Performing a type 2 deformation pushes part of the loop into the regular neighborhood.
Keeping in mind that the aim of a type 2 deformation is to kill off an intersection with a particular surface, it is not hard to see that we should push the arc through the other side of $N(F_1)$.   Note that since we will only carry out such moves when we find a type 1 reduction, the following algorithm will still be finite time.
We will examine these moves for  the checker-board hierarchies in more detail in section~\ref{sec:moves}.

\medskip

The following algorithm determines whether a loop $l \subset S^3 - N(L)$ is contractible and hence solves the word problem:

\begin{Alg}[\cite{Du1982}]  \label{surfalg}
If $l \cap F_1 = \emptyset$ then use process X to check for type 1 reductions with respect to $F_2, \ldots , F_n$ and carry them out (since $l \cap F_1 = \emptyset$ there are no possible type 2 deformations).  If now  $l_j \cap \cup_{i=2}^n F_i = \emptyset$ then $l$ is contractible, otherwise it is not.

 If $l \cap F_1 \neq \emptyset$ then use $F_1$ to split $l$ into arcs $l_1 , \ldots , l_m$. Use algorithm~\ref{f1alg} to check for and carry out all type 1 reductions with respect to $F_1$ for each arc $l_i$, $i=1, \ldots m$.
If now $l \cap F_1 \neq \emptyset$ then $l$ is not contractible.

Otherwise split $l$ into arcs $l_1, \ldots ,l_s$ by the intersections with $F_2, \ldots , F_n$ and use process X to check for, and then carry out, all type 1 reductions.
If now  $l \cap \cup_{i=1}^n F_i = \emptyset$ then $l$ is contractible, otherwise it is not.
\end{Alg}

\sgap

\section{Loops and the Dehn Presentation} \label{sec:loo}

We will now  look at the relationship between the (augmented) Dehn presentation and loops in the link complement. This provides the interaction between small cancellation theory and normal surfaces.

By a loop we will mean either the path itself or its embedding. It will be clear from context which is meant, and since we will mostly work up to homotopy, this confusion does not cause any problems.

Let $l$ be a loop in $S^3 -L$.
   Up to homotopy we may assume that $l$ intersects the interiors of the  regions of the link transversally. We will always assume this of any loop or arc.
 Further, since the union of the regions splits $S^3$ into two 3-balls, every  arc between two intersection points whose interior does not intersect any regions, determines a unique homotopy class.
So, up to homotopy, the only information that a particular (oriented) loop carries is which regions, in what order and in what direction it intersects these regions.

Of course we can consider a word in the (augmented) Dehn presentation as a based oriented loop in the (augmented) link complement.  What about the other direction?

If an oriented  loop is based then, it is equivalent to a word in the alphabet generated by the labels of the regions, $x_0, \ldots , x_n$, constructed by
following the loop in the direction of the orientation from the base point and assigning the letter $x_i$ every time the loop passes downward through the region $x_i$  and $x_i^{-1}$ every time the loop passes upward through the region.
We call the word generated in this way the {\em canonical word} determined by the loop.
Note that all canonical words  are even in length and alternating in sign.

We have constructed a 1-1 correspondence between alternating words of even length
 and the homotopy classes of based oriented loops in $S^3 -L$ relative to the set of intersection points with the union of the regions of the link.
Let's see how this relates to the Dehn presentation.

Given a based oriented  loop $l \subset S^3-L$, we would like to obtain a canonical element in the augmented Dehn presentation.
The usual way to do this is to pull everything upwards, by homotopy, to the base point. This way is not suitable for our purposes as it introduces extra intersection points with the regions  of the augmented link.
Instead we prefer the following method.
We replace each arc from a region $x_i$ to a region $x_j$ which lies completely above the regions with  a path from the intersection point with $x_i$ directly to the base point of the space and then directly back down to the intersection point with $x_j$.
Each arc defined by the intersection of the loop with the regions which lies entirely underneath the regions, we homotope this by sweeping it around the outside of the $\mathcal{O}$ component of the augmented link so that it intersects the base point to the space at one point.
Up to homotopy we may assume that the base point of the loop coincides with that of the space.
It is clear how this construction represents an element in the Dehn generators of the augmented link group.
We call this the {\em canonical element} of the augmented Dehn presentation determined by $l$.
It is easy to see that the canonical element of a loop is exactly the element in the augmented Dehn presentation given by the canonical word.

Of course this is not well defined.  The problem being that if  we are given two freely homotopic based oriented loops in the link complement, taking the canonical element fixes a path from the loop to the base point of the space and as elements of the fundamental group  and these may represent different elements of the fundamental group.  However, it is not hard to see that these two classes are conjugate.  Therefore we have a correspondence between based oriented loops in the link complement and included words in the augmented Dehn presentation which is well defined up to conjugacy.  Since we are interested in the word and conjugacy problems, this is sufficient for our purposes.

\begin{Remark}
This is basically the folklore result that, for sufficiently nice topological spaces, the conjugacy problem  is equivalent to determining whether two loops are freely homotopic.
\end{Remark}

\section{Hierarchies and the Dehn Presentation} \label{sec:moves}
In this section  $L$ is an alternating link.  We will explore the relationship between the checker-board hierarchy and the augmented Dehn presentation.

We constructed the black (resp. white) checker-board hierarchy from the regions of the link by identifying the black (resp. white) regions and removing a neighborhood of the boundaries of these surfaces.  So we have a correspondence between the hierarchies and the regions. Furthermore, we may assume that the intersection points between the loops and the regions lie away from the parts of the regions which are removed during the construction, giving a correspondence between loops in the  link complement containing a black or white checker-board hierarchy  and the set of regions,  which, by section~\ref{sec:loo}, gives a correspondence with words in  the augmented Dehn presentation.
From this, one expects  a correspondence between type 1 and type 2 moves and some actions on the augmented Dehn presentation.  We will work out the details.

Although a type 1 reduction is performed with respect to a surface in the hierarchy, rather than the regions of the link, the two end points of the arc we are reducing must lie in the same region otherwise that arc would intersect another surface in the checker-board hierarchy.
It is now easy to see that in terms of regions, a type 1 reduction  is a move which pulls an arc which intersect a region $x_i$ then comes straight back through in the opposite direction, completely through the region. In terms of the canonical words this corresponds to replacing $x_i^{\pm 1}x_i^{\mp 1}$ with $1$, which is a free reduction.
Similarly, a type 1 augmentation with respect to $x_i$ introduces a subword $x_i^{\pm 1}x_i^{\mp 1}$.

\begin{figure}
\centering
\epsfig{file=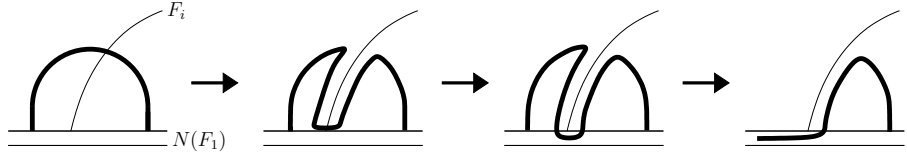, width=12cm}
\caption{A type 2 deformation.}
\label{finger2}
\end{figure}

\begin{figure}
\centering
\epsfig{file=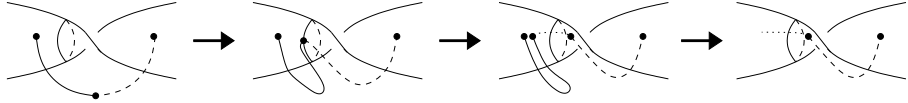, width=12cm}
\caption{Realizing a type 2 deformation.}
\label{celltype2}
\end{figure}

\begin{figure}
\centering
\epsfig{file=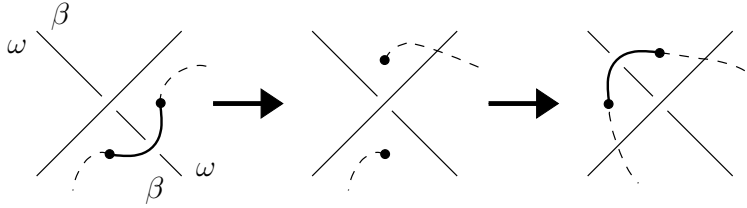, width=10cm}
\caption{A projection of a type 2 deformation.}
\label{type2image}
\end{figure}

Let's look at the more complicated type 2 deformation.
 A type 2 deformation can be thought of as a finger move which pushes an arc along $F_i$, $i \geq 2$, into $ N(F_1)$, then push everything on one side of $F_i$ into the neighborhood as in figure~\ref{finger2}.
Remembering that the hierarchies are formed by regions which intersect as in figure~\ref{2-cells}, we can interpret this in the link complement as
 a finger move which pushes the intersection point along $F_i$ until it intersects  the North-South axis. Then pushing it into $N(F_1)$, so that it does not intersect with $F_i$ at all.   Finally, push the relevant part of the arc so that part lies completely inside the neighbourhood. This is indicated in  figure~\ref{celltype2}, where the dotted line may intersect more  surfaces before returning to $F_1$.
Next we go on to push the arc lying in the neighbourhood through to the other side.  This process is shown as a projection in figure~\ref{type2image} (where the type 2 deformation is with respect the black checker-board hierarchy).
By the correspondence with the regions we see that a type 2 deformation of some loop corresponds to exchanging a pair in the canonical word.

Our final observation is that moving the base point of a based oriented loop just changes where we start (and finish) reading the canonical word and therefore corresponds to a cyclic permutation of the word.

\section{A Comparison of the Algorithms for the Word Problem} \label{sec:comp}

We will apply the theory of sections~\ref{sec:loo} and \ref{sec:moves} to the two solutions of the word problem for prime alternating link groups.
We will begin by expressing   Dugopolski's algorithm in terms of  exchanging pairs and free cancellations of words in the augmented Dehn presentation.

In short, Dugopolski's algorithm takes a loop in $S^3-L$ and checks for and carries out type 1 and 2 moves to determine whether it contracts or not.  We  express this procedure in terms of words in the augmented Dehn presentation.

\smallskip

Recall that in algorithm~\ref{f1alg} we used process X to check for and carry out a consecutive sequence of type 2 deformations which may result in a type 1 reduction.  This process has the following interpretation:

Given a word $w = v(1)w(1)w(2)\cdots w(n)v(2)$, such that $v(i)$, $i=1,2$ are of parity black and $w(i)$, $i=1,\ldots , n$ are of parity white, define {\em process Y} according to the following method: if $v(1)w(1)$ is a pair then it uniquely determines a relator in the presentation, replace this pair with the remaining letter of parity black in the relator.  Repeat this for the first two letters of the new word  obtained. Repeat this process until we reach a word where the first two letters are not a pair.

If process Y returns a word of two letters (both of parity black) which freely cancel we say that it was {\em successful}.

If $w$ is a word for which process Y was successful we define the {\em switch} of $w$ to be the word obtained as in process Y but by substituting each pair with the other pair from the relator it determines and canceling the last two letters once the process has been completed.
A switch corresponds to carrying out all the sequence of type 1 and 2 moves found by algorithm~\ref{f1alg}.

\begin{Lem}
A switch is a chain collapse.
\end{Lem}
\begin{proof}
If process Y is successful we can form a chain whose chain word is $v(1)w(1)w(2) \cdots w(n)v(2)$.  Each substitution in the formation of the switch replaces the black-white path with the white-black path around each relator square.  Finally the free cancellation kills off the remaining `spike'.
\end{proof}

Given a word $w = v(1)w(1)w(2)\cdots w(n)v(2)$, such that $v(i)$, $i=1,2$ are of parity black and $w(i)$, $i=1,\ldots , n$ are of parity white, we can rewrite algorithm~\ref{f1alg} as:

\begin{Alg}\label{type1words}
If $n=0$ check whether $v(1) \eqpi v(2)^{-1}$.  If $n \neq 0$ carry out all free reductions on the subword $w(1)w(2)\cdots w(n)$.  If now $n=0$ check if $v(1) \eqpi v(2)^{-1}$. If $n \neq 0$ do process Y.

If we were successful at any point then a type 1 reduction with respect to $F_1$ is possible, otherwise one is not.
\end{Alg}

\sgap

Now let $l$ be an element of $\pi_1 (S^3 -L)$, where $L$ is an alternating link. Without loss of generality, we may assume that  $\pi_1 (S^3 -L)$ is given by the Dehn presentation.  Then let $w$ be the included canonical word representing this loop.

The following is a rewriting of algorithm~\ref{surfalg} and hence it solves the word problem for the groups of alternating links.

\begin{Alg}\label{dugoword}
If $w$ contains letters of parity white only then freely reduce $w$.  If this gives the empty word then $w \eqpi 1$ otherwise $w \neqpi 1$.

If $w$ contains letters of parity black. Then if there is  exactly one letter of parity black, $w \neqpi 1$, otherwise cyclically permute so that the first letter is black.  Now $w$ can be written as $\beta_1 \omega_1 \beta_2 \omega_2 \cdots \beta_n \omega_n$, where $\omega_i$ is a subword of white parity (and may be empty), and each $\beta_i$ consists of exactly one letter of parity black. Pick a subword $\beta_i \omega_i \beta_{i+1}$, for $i$ modulo $n$, and run algorithm~\ref{type1words}.  If the algorithm is successful freely reduce the subword or replace it with its switch as appropriate.  Repeat this with the new word until no more such substitutions are possible.

If now there are only parity black letters then $w \neqpi 1$.  If there are only parity white letters then freely reduce.  If the word is now empty then $w \eqpi 1$ otherwise $w \neqpi 1$.
\end{Alg}

\sgap

This simplifies to the following:
\begin{Duwo}\label{dugowordsmall}
Freely reduce $w$ and its cyclic permutations. Chain collapse all chain words  beginning  and ending in a letter of parity black in $w$ and its cyclic permutations.  Freely reduce.  If we obtain the empty word then $w \eqpi 1$ otherwise it does not.
\end{Duwo}

\sgap

Observing that we need never cyclically permute the word, we see that this is exactly algorithm~\ref{dugoword}, where the chain collapses are only carried out with respect to black chains.

\medskip

We now turn our attention to the rewriting of the group theory algorithm  in the language of normal surfaces.

We start with some included word $w$ in the augmented Dehn presentation.  Algorithm~\ref{smallalg} tells us how to determine whether this represents the identity in terms of free cancellation and chain collapses. We know that a free reduction corresponds to a type 1 reduction. Let's see what we can do with a chain collapse.

Let $a t b$ be a chain word  with inner link path $s$, where $t=t_1 t_2 \cdots t_n$ and   $s=s_1 s_2 \cdots s_n$.
Then we may view a chain collapse as a sequence of exchanging pairs followed by  a free reduction:
$ a t_1 \cdots t_n b \rightarrow s_1 u_1 \cdots t_n b \rightarrow \cdots \rightarrow
s_1 \cdots s_n b^{-1} b \rightarrow s_1 \cdots s_n$, where the $u_i$ are the labels of the internal vertical edges in the chain.
We want to write this as a sequence of moves on a loop. Recall that the hierarchies split the loop into subarcs.  Starting from the intersection of the loop with region $a$ do a type 2 deformation (with respect to the checker-board hierarchy with the same colour as $a$) on the subarc between regions $a$ and $t_1$. Pushing this through the surface our subarc now intersects regions $s_1$ and $u_1$.  Now do the type 2 deformation on the next  subarc between $u_1$ and $t_2$.  Continue like this until we reach a subarc between the region $b$ (this is after $n$ moves) and eliminate this through a type 1 reduction.

Letting $F_1 , \ldots , F_n$ be the black checker-board hierarchy and  $H_1 , \ldots , H_n$ be the white checker-board hierarchy we can interpret algorithm~\ref{smallalg}  as:

\begin{Alg} \label{wordhiealg}
If $l \cap \cup_{i=1}^n F_i = \emptyset $ then $l$ is contractible.  If $l \cap \cup_{i=1}^n F_i \neq \emptyset $ use process X to find and carry out all type 1 reductions except for those which would move the base point. If now $l \cap \cup_{i=1}^n F_i = \emptyset $ then $l$ contracts. Otherwise use the black hierarchy to split $l$ into arcs $f_1, \ldots , f_p$ and  use the white hierarchy to split $l$ into arcs $h_1, \ldots , h_q$. Starting from the base point use algorithm~\ref{f1alg} to search for a type 1 reduction with respect to both $F_1$ and $H_1$. If any are found which do not include the base point of $l$, look at the arcs of $l$ which they determine.  These will be nested, so choose the innermost arc and carry out the sequence of moves required to make the relevant type 1 reduction.  Use process X to find and eliminate any further type 1 reductions which do not move the base point. Repeat this until all such reductions are made.  If now $l \cap \cup_{i=1}^n F_i = \emptyset $ then $l$ contracts.  Otherwise it does not.
\end{Alg}

\sgap

Note that in this algorithm, we keep the base point fixed at all times and we are working with respect to two hierarchies.

\begin{Remark}
By the geodesic characterization theorem, algorithm~\ref{wordhiealg} will  find a loop in the homotopy class of $l$ such that the number of intersection points with the regions (or hierarchy) is minimal.
\end{Remark}

\section{The Role of Non-Positive Curvature}

The role of non-positive curvature  in small cancellation theory is well known. In this section we will discuss where the non-positive curvature is hidden in the normal surface approach. This suggests why the two algorithms coincide on prime alternating links.
The 2-complex discussed here was noted independently by Aitchison (unpublished) and Wise in \cite{Withesis} (see also \cite{Brbook}).

Given a link $L$  construct a PE 2-complex as follows: take two 0-cells, $v_-$ and $v_+$. Add a 1-cell corresponding to each region of the link, oriented from $v_+$ to $v_-$ and label each 1-cell according to the regions. Finally take one 2-cell $D_i$ for each relator $r_i$ and attach the 2-cells by a continuous map taking $\partial D_i$ to the path in the 1-skeleton representing $r_i$.
Give this complex a PE structure by regarding the 2-cells as regular 4-gons.
We call this the {\em Dehn complex} of $L$ and we call the Dehn complex of the augmented link the {\em augmented Dehn complex} of $L$.

Clearly the fundamental group of the (augmented) Dehn Complex of a link $L$ is the (augmented) Dehn presentation of $L$.

\begin{Lem}
Let $L$ be a reduced link. Then the Dehn complex embeds into $S^3-L$.
\end{Lem}

\begin{proof}
Recall that we are assuming $L \subset \mathbb{R}^3 \cup \infty$ coincides with its projection except at crossing balls.  Take $v_+$ to lie above $\mathbb{R}^2 \cup \infty$ and $v_-$ to lie below.  Position the 1-cells so that the pass through the region which label them. The result follows from the following picture which shows the embedding locally at a crossing.
\[
\epsfig{file=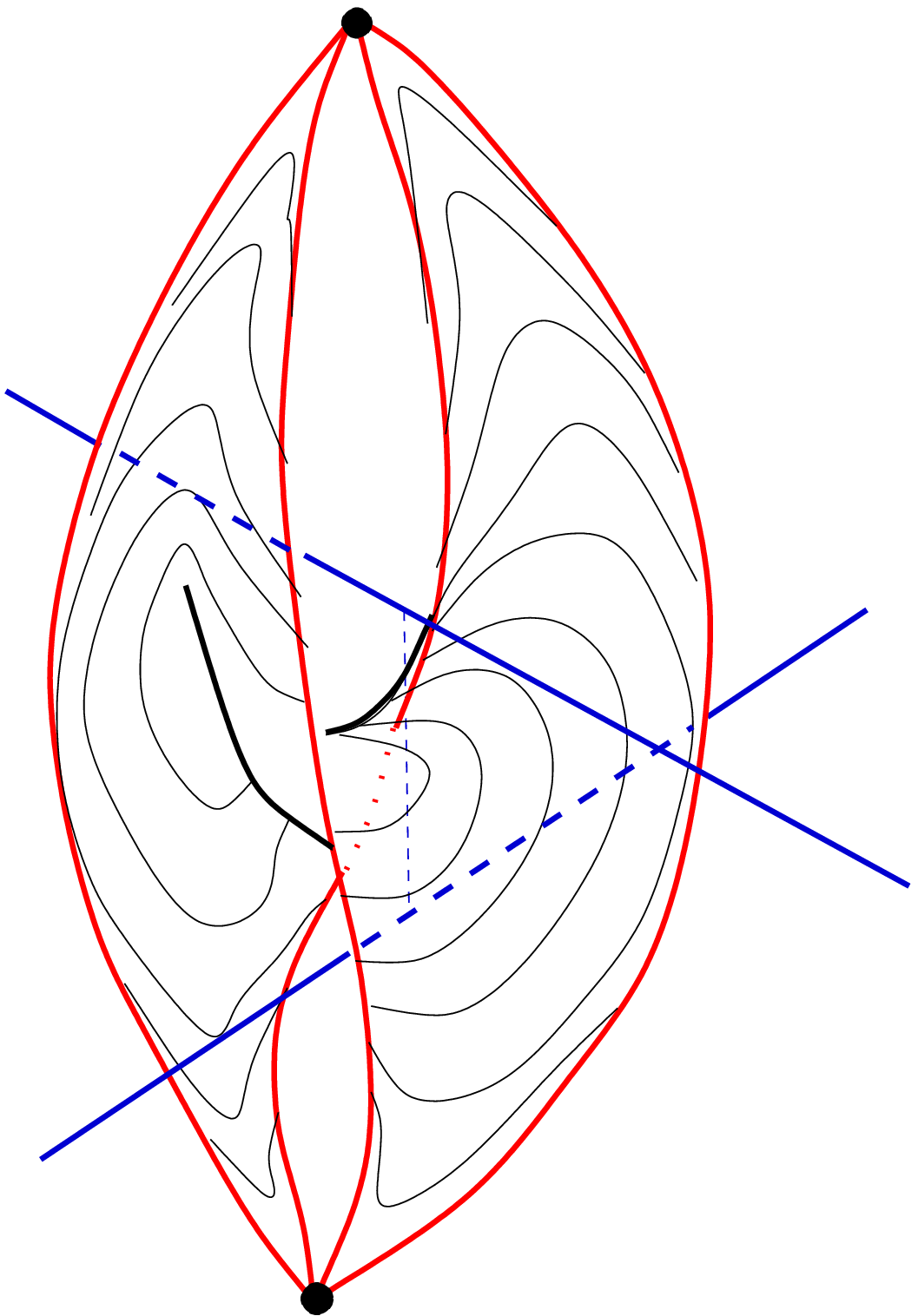, height=6cm}
\]
\end{proof}

\begin{Remark}
With a little more thought one can see that when $L$ is reduced the Dehn complex is in fact a deformation retract of the link complement.
\end{Remark}

Consider the Dehn complex embedded as in the above lemma.
Recall the discussion of section~\ref{sec:loo} which says that up to homotopy the relevant  information of any loop are its intersections with the regions.  So canonically we may assume that any loop coincides with the 1-skeleton of the Dehn complex.  We may think of type 1 reductions and type 2 deformations as moves between paths in this 1-skeleton.  In this situation, a type 1 reduction is clearly a move within the one skeleton.  A type 2 deformation corresponds to exchanging a pair and it is easy to see that we may assume that this deformation lies on the appropriate 2-cell of the Dehn complex.
Thus we see that Dugopolski's algorithm descends to an algorithm on the Dehn complex of a reduced alternating link.

The small cancellation properties of the \adp are dependent upon the Dehn complex.
We say that a PE 2-complex has {\em non-positive curvature} if traveling around the boundary of any small disk in the complex whose origin is at a 0-cell we travel through an angle of $\geq 2\pi$.

The following proposition is theorem~\ref{dehncan} rephrased in terms of non-positively curved complexes.
\begin{Prop}[\cite{{We1971},{Withesis}}] \label{nonposcurv}
The Dehn complex of a reduced link is non-positively curved if and only if the projection is prime and alternating.\end{Prop}

One can make any non-positively curved square complex into one with a \c4t4 fundamental group by adding a 1-cell between a distinguished 0-cell and every other 0-cell in the complex.
These observation and the above proposition indicate why the two approaches discussed in this chapter coincide.

\section{A Small Cancellation Solution for the Conjugacy problem} \label{sec:smallconj}
In this section we will outline Johnsgard's solution of the conjugacy problem.  Johnsgard analyzes the types of disc diagrams that are possible for \c4t4 small cancellation groups and the restrictions imposed upon the diagrams by the geometry of the link.  An algorithm is then provided for finding  cyclic geodesic representatives of the conjugacy class of a given element.  The reader is referred to \cite{Jo1997} for the details and justification of the results discussed in this section.

\sgap

First some nomenclature for  \c4t4 presentations.  If there exists a letter $b$ such that $ab$ and $b^{-1}c$, where $c \neq a^{-1}$, are both pairs, we call $ac$ a {\em sister-set}.   If a disk of relator squares forms a ``L'' shape such that the inner crook of the ``L'' is not a pair, we call the word labeling the crook a {\em pseudo-pair}.
Note that in a \cant4 presentation no pair is a sister-set  and in a group presentation with parity, no pseudo-pair is a sister-set.

We say that a word is {\em cyclically geodesic} if all its cyclic permutations are geodesic.  If a word is a cyclic permutation of another, we say that the two words are {\em cyclically equal},
and if two geodesic words represent the same group element we say they are {\em equivalent geodesics}.

\sgap

We will use the following construction in algorithm~\ref{joconjalg}. Regard the integral points $(n,m)$ in the plane as potential 0-cells for a 1-complex.
Let $w$ be a geodesic word in a \c4t4 presentation.  We embed $w$ in the plane in the following way.  Start at the point $(0,0)$, which we take to be a o-cell. Choose a point at displacement $\vec{i}$ away and also take the 1-cell between these points as part of the complex.
 Label this 1-cell with the first letter of $w$.
  If $w$ has more letters then from this point choose a 0-cell and corresponding 1-cell a displacement of $\vec{i}$ away and label this new 1-cell with the next letter in the word according to the following criterion: if the new letter and the preceding one form a pair or a pseudo-pair which is not a sister set, change direction from that most recently taken; otherwise keep going the same way.  Continue in this way for the remainder of the word.  The 1-complex thus obtained is called the {\em standard embedding} of $w$.

Denote the standard embedding by $C_0$.  From $C_0$ we construct a singular disc diagram  by the following iterative procedure: consider each right angle in the plane that is not contained in consecutive edges of a relator square of $C_i$, if this right angle is labelled by a pair then add the (unique) relator square to the diagram it determines.  Call the resulting diagram $C_{i+1}$.  Continue in this way until no more squares can be added.

 The resulting diagram is called the {\em geodesic completion} of $w$.  It is characterized by the following theorem:

\begin{Geocomthe}[\cite{Jo2000}]
A geodesic word in a \c4t4 presentation uniquely determines a square tiling (the geodesic completion) bounded by  a rectangle in the Euclidean plane such that every equivalent geodesic of the original word is a label of a geodesic rectilinear edge path which is path-homotopic to the edge path of the original word.
\end{Geocomthe}

The geodesic characterization theorem provides an easy method for replacing a word with a conjugate cyclic geodesic word by taking cyclic permutations and geodesic representatives of the word in \c4t4 presentations.

We also note that by \cite{Jo1997} proposition 4.1 and corollary 4.8, in an alternating \c4t4 presentation  with parity,  conjugate cyclic geodesic words are of the equal even length, both alternate in sign, and have equal numbers of letters of each parity.

\sgap

Let $u^{\prime}$ and $w^{\prime}$ be words in the augmented Dehn presentation of a prime alternating link. Then the following algorithm solves the conjugacy problem in polynomial (order 7) time.

\begin{Alg}[\cite{Jo1997}] \label{joconjalg}
Take cyclic geodesic representatives $u$ and $w$ of $u^{\prime}$ and $w^{\prime}$. We may assume these are non-empty words.

If $u$ and $w$ are not both of equal even length,  alternating in sign with the same number of letters of each parity, then $u^{\prime}$ and $w^{\prime}$ are not conjugate.
Otherwise there are two cases: whether $w$ has letters of both parities or not.

First consider the case where $w$ has letters of both parities.
We are going to construct a planar 2-complex using  relator squares of unit length coming from the presentation and regarding the integral points of the plane as possible 0-cells.

Begin by taking the geodesic completion of $w$.  For simplicity, assume that the first letter of $w$ is oriented in the $\vec{i}$ direction and the path representing $w$ lies in the first quadrant with endpoint $(I,J)$.

Some of the relator squares in the geodesic completion may intersect the lines $x=0$ or $x=I$.  Let $l_1$ (resp. $r_1$) be the word labelling the path in the geodesic completion which sits on the line $x=0$ (resp. $x=I$).
If $l_1$ is non-empty then add a path from $(I,J)$ to $(I,J+|l_1|)$ labelled by $l_1$ to the complex and take the geodesic completion.
If this creates a new path in the geodesic completion which lies on $x=0$ (so it initial point is $(0,|l_1|)$),  we label this $l_2$.  We place a copy of this path on top of $(I,J+|l_1|)$ and take the geodesic completion.  Continue in this way until we either reach an $l_i$ which is empty or until the 2-complex starts repeating itself.

Repeat this process with the words $r_i$ which we add to the complex on the $- \vec{j}$ direction from the point $(0,0)$.

Now check every path from point $(0,n)$ to $(I, J+n)$, where $n$ is an integral point on the 2-complex, for cyclic permutations of $u$.  If such a path is found $u$ and $w$ are conjugate otherwise they are not.

Now consider the case where  $w$ is composed entirely of letters of a single parity.
Look for all letters $x_{a_i}$ such that $x_{a_i}^{\pm 1} \cdot w \cdot x_{a_i}^{  \mp 1}$ is a chain word.
This determines some chains which all have a side labelled by $w$.  Identify these chains along their common side to obtain a connected 2-complex.
This 2-complex has a set of paths which are the sides of some chain.
For each of these paths we look for all chains which have a side labelled by a path in this set and, unless it is
 a cyclic permutation of some chain which has been previously added to the complex, we identify the common edges of the chain and the 2-complex.  Continue in this way until no new chains can be added to the complex.

If there is a path in this complex which was the side of a chain which is labelled by a cyclic permutation of $u$ then $u$ and $w$ are conjugate otherwise they are not.
\end{Alg}

\sgap

\begin{Remark}
We will see in section~\ref{app:planar} that the 2-complexes produced by the algorithm are in fact simply connected.
\end{Remark}

\section{The Algorithm in Action}
To aid the digestion of the conjugacy algorithm, we provide a few examples of the diagrams it produces.

First consider the trefoil.  This has augmented Dehn presentation
\[
\langle x_0, \ldots , x_4 | x_1 x_4^{-1} x_2 x_0^{-1} ,
 x_1 x_0^{-1} x_3 x_4^{-1} ,
x_2 x_4^{-1} x_3 x_0^{-1}
\rangle.
\]
The following diagrams are constructed by the algorithm for  the words $w= x_3^{-1} x_4 x_0^{-1} x_2$ and $v=x_1 x_3^{-1}$.
The diagram for $w$ is on the right and for $v$ the left.
\begin{center}
\epsfig{file=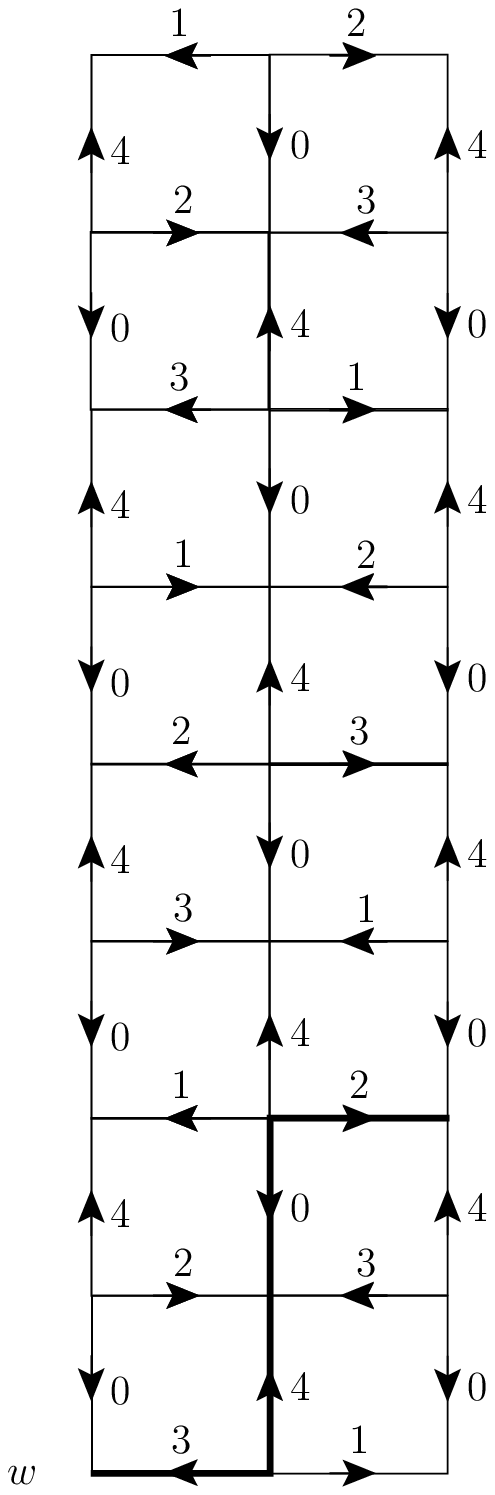, height=9.6cm}
\hspace{3cm}
\raisebox{1.2cm}{\epsfig{file=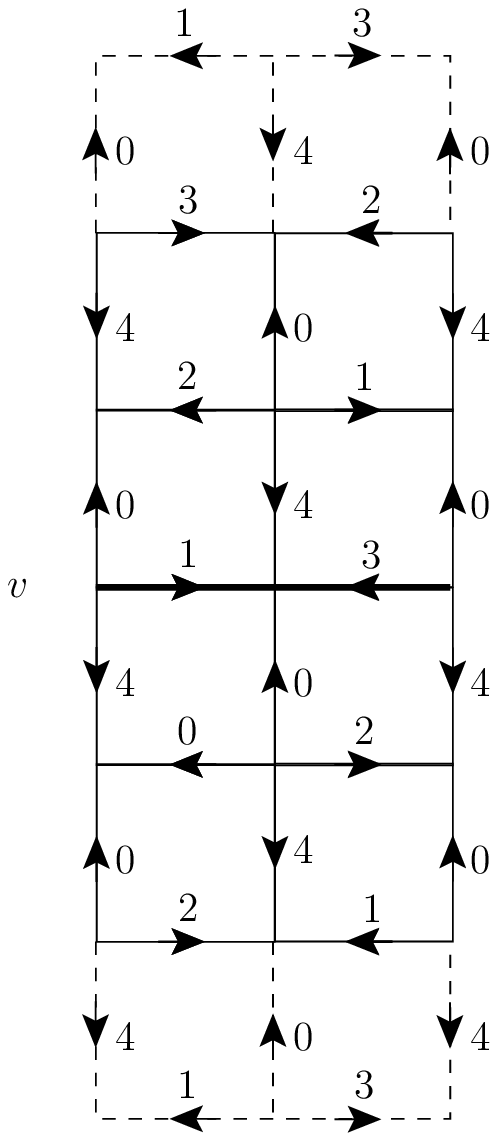, height=7.2cm}}
\end{center}
Reading off the diagrams we see  that $w$ and, say, $ x_2^{-1} x_3 x_0^{-1} x_4$ are conjugate and that any two words with these as their cyclic geodesic representatives are conjugate.  It also proves that say, $ x_1^{-1} x_4 x_2^{-1} x_4$ is not conjugate to $w$.
Similarly, the right diagram tells us that $u$ is conjugate to $ x_2^{-1} x_3 $ and is not conjugate to
$ x_4^{-1} x_0 $.

\sgap
The tiling for $w$ above fills a  rectangle.  In general this wont happen.  To see an example of how pseudo-pairs can change things, consider the figure eight knot (the augmented Dehn presentation of the trefoil has no pseudo-pairs).
Its augmented Dehn presentation is:
\[
\langle x_0, \ldots , x_{5} | x_1 x_4^{-1} x_2 x_0^{-1} ,
 x_1 x_0^{-1} x_3 x_4^{-1} ,
x_2 x_4^{-1} x_3 x_5^{-1},
 x_2 x_5^{-1} x_3 x_0^{-1}
\rangle.
\]
We construct the  diagram for $w= x_2^{-1} x_1 x_{5}^{-1} x_3$, noting that $x_1 x_{5}^{-1}$ is a pseudo-pair.
\begin{center}
\epsfig{file=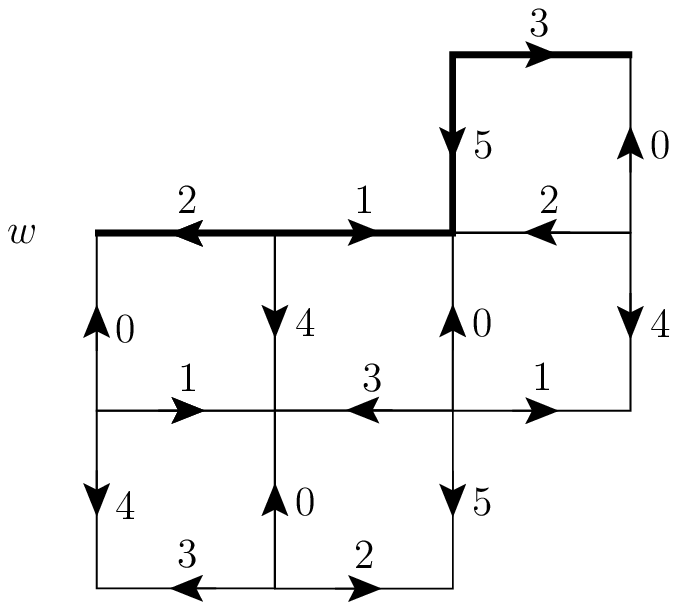, height=4cm}
\end{center}

\sgap

For further example of how things can differ from the above, consider the link below (whose lengthy presentation we exclude).
\begin{center}
\epsfig{file=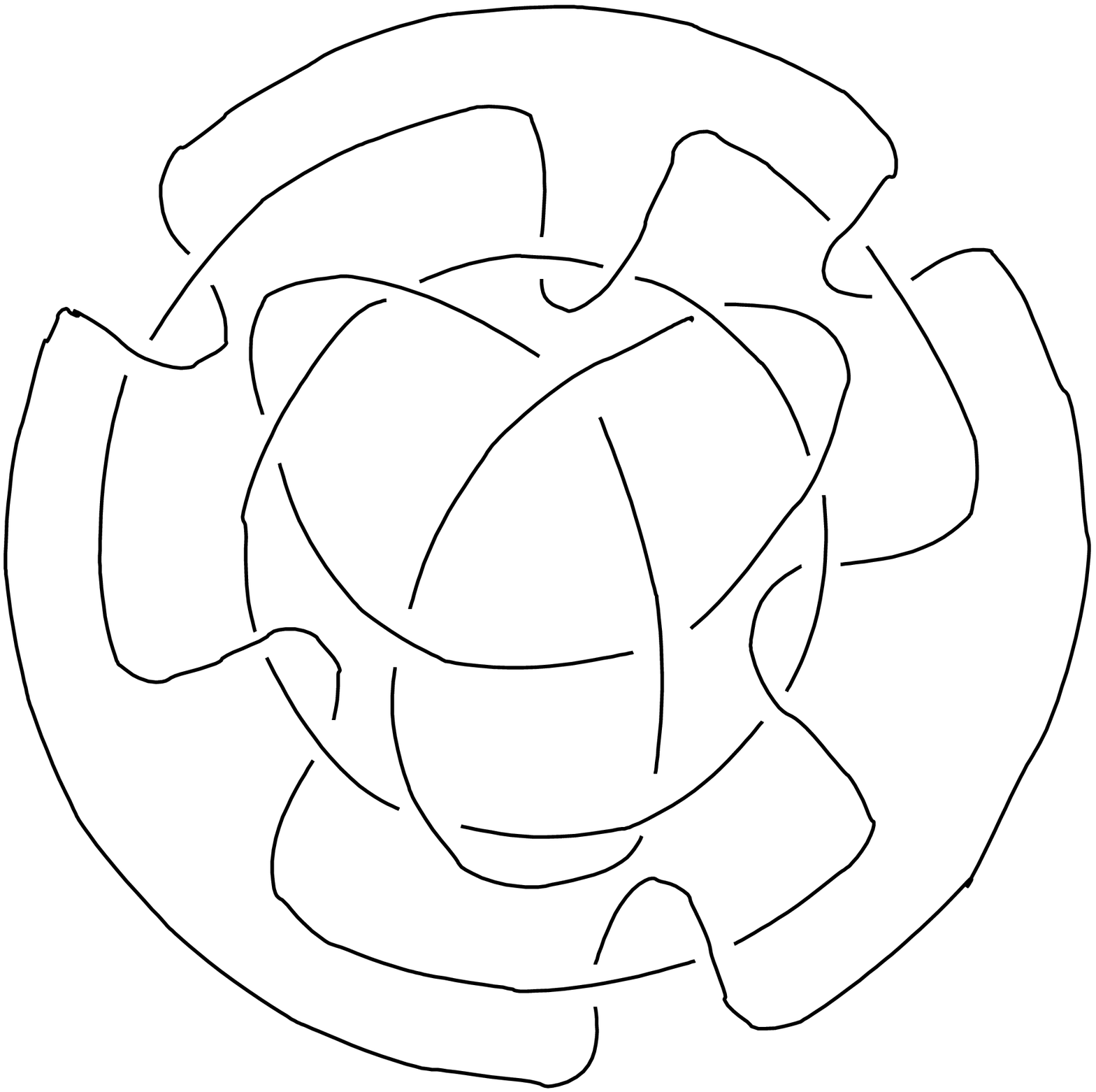, height=3cm}
\hspace{2cm}
\epsfig{file=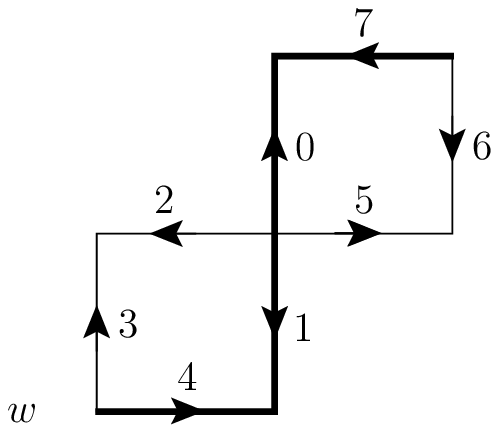, height=3cm}
\end{center}

There exists a labelling of the regions of the above link such that $x_2 x_1^{-1} x_4 x_3^{-1}$ and $x_7 x_6^{-1} x_5 x_0^{-1}$ are relators and $ x_2^{-1} x_0$,  $ x_1^{-1} x_5$, $ x_7^{-1} x_3$ and $ x_6^{-1} x_4$ are not pairs or pseudo-pairs.  In fact in this labelling  1 and 0 are of opposite parity.

Notice that the diagram produced by algorithm~\ref{joconjalg} from the word $w=x_4   x_1^{-1} x_0   x_7^{-1}$ contains a split vertex and the diagram contains collinear edges with opposite parity.
In fact, by the theory of \cite{Jo1997} in such situations (where there is a two letter subword which is not a pair, sister set or pseudo-pair) all of the cyclic geodesic representative come from the geodesic completion.

\section{Normal Surfaces and the Conjugacy Problem} \label{sec:normcon}

In light of what has come before, it should come as  no surprise that the small cancellation solution to the conjugacy problem gives a normal surfaces algorithm.
Algorithm~\ref{joconjalg} finds finite sets of certain cyclic geodesic representatives of the conjugacy classes of the group. Moreover, every group element is conjugate to one of these representatives.

The algorithm  starts with an included word $w= i_* (w^{\prime})$, where $i_*$ is the inclusion of the Dehn presentation into the augmented Dehn presentation, and constructs a planar diagram. From this diagram we can read off a set representatives of the conjugacy classes of $w$.
By section~\ref{sec:loo}, $w$ can be represented as a based oriented loop in $S^3-L$ and since we are interested in conjugacy, by section~\ref{sec:comp}, we can work with these loops.
In this section we will interpret the small cancellation methods for the conjugacy problem to find a normal surfaces solution to the conjugacy problem for prime alternating links.  Our aim is to construct the (representatives of the) conjugacy classes geometrically.

\sgap

 We call a loop in $S^3 -L$ or $S^3 -(L \cup \mathcal{O})$ {\em minimal} if it attains the minimal number of intersections with the regions $x_0 , \ldots , x_n$ (or equivalently, for $S^3-L$, the hierarchy).
Say that two loops are {\em equivalent minimal} if they are homotopic and minimal.

It is clear that  canonical loops from  cyclic geodesic words are minimal and loops of equivalent geodesics, cyclic permutations and inverses from a cyclic geodesic word are equivalent minimal.

\sgap

First we need to  understand the geometry behind the geodesic completion.  Remember that by the geodesic completion theorem,  the geodesic completion of a word $w$ contains all of its equivalent geodesics.  One point to bear in mind, is that different cyclic permutations of the same word will in general produce different geodesic completions (unless there are no pseudo-pairs).  This means that our interpretation to loops will depend upon the base point.

By the construction of the geodesic completion, it is clear that the following method will produce the set of all equivalent geodesics of a word $w$: Start with the set containing only $w$.  Scan $w$ for pairs. To each pair we find, take a copy of $w$ and replace that pair with the other pair in the relator it determines. If this creates a word not already in the set add it to the set.  Repeat this process with each new word we find which is not already in the set.  Stop when we can't continue with this process.

\sgap

Note that if $w$ has $n_{\omega}$ of letters of parity white and  $n_{\beta}$ of parity black, then there are at most $(n_{\beta} + n_{\omega})! / n_{\beta}! n_{\omega}!$ equivalent geodesics.

Now since exchanging a pair corresponds to a type 2 deformation of the canonical loop, it is easily seen from the construction that the set of (minimal) loops produced is the set of loops which can be obtained from $l_w$ by a sequence of type 2 deformations which keeps the base point fixed. Formally:
\begin{Lem}
Let $l_{i_*(w)}$ be the oriented based loop representing  a geodesic included word $i_*(w)$  in the augmented Dehn presentation of a prime alternating link.
Construct a set of based oriented loops by performing all possible sequences of type 2 deformations while keeping the base point fixed. This set represents the equivalent geodesic words of $w$.
\end{Lem}

We move on to the question of generating the  set of representatives found by the algorithm in terms of geometric moves on loops.  Denote by $l_w$ the based oriented loop determined by the word $w$ from the inclusion of $\pi_1 (S^3 -L)$ into the augmented Dehn presentation.

The first step of algorithm~\ref{joconjalg} is to find cyclic geodesic representatives of the words.  Algorithm~\ref{wordhiealg} will do this.

The next step of the algorithm concerns words which are of odd length or non-alternating. Since words coming from the inclusion of the Dehn presentation into the augmented Dehn presentation are necessarily alternating and of even length, we need not consider this case.

Now suppose that $w$ is alternating, of even length and contains letters of one parity only.  Then algorithm~\ref{joconjalg} produces a finite strip in the plane. The set of representatives are the sides of the chains in this strip and the set of all geodesic representatives of the conjugacy class are their cyclic permutations.  The set of representatives produced by the algorithm determines a set of minimal based loops.  We would like to understand how to generate these geometrically starting from $l_w$.  To do this it is sufficient to understand how to move between the sides of a single chain of the construction.

Consider a conjugacy chain with chain word $a t a^{-1}$ and inner link path $s$, where $t=t_1 t_2 \cdots t_n$ and   $s=s_1 s_2 \cdots s_n$.  So the sides $s$ and $t$  determine loops $l_s$ and $l_t$ respectively and we want to obtain $l_t$ from $l_s$.  On the level of groups, we can make the  sequence of substitutions:
 $ s_1\cdots s_n \rightarrow a^{-1} s_1 \cdots s_n a \rightarrow  t_1 b_1 s_2 \cdots s_n a \rightarrow  t_1 t_2 b_2 s_3\cdots s_n a \rightarrow \cdots \rightarrow   t_1 \cdots t_n a^{-1} a \rightarrow t_1\cdots t_n $.
From this we get the following geometric interpretation:
Deform the arc of the link containing the base point so that the base point lies over region $a$
 and no further intersections with regions are added.
  Further deform this arc by a type 1 augmentation with respect to $a$ (ie push the base point through region $a$).
Follow the loop from the base point in the direction of the orientation. Each time we meet an intersection point carry out a type 2 deformation on the arc between that intersection point and the following one and continue following the loop.  Continue in this way until we have carried out $n$ deformations (or equivalently we have reached the final intersection point before the base point). Now carry out the type 1 reduction with respect to region $a$ (since our original loop was minimal there is only one choice for the reduction).  The reductions are with respect to the black or white checker-board hierarchy depending upon the parity of $a$.

Note that two loops are related by the above sequence of moves if and only if they can be represented as two sides of a chain.

So by algorithm~\ref{joconjalg} we have:
\begin{Lem} \label{normsingle}
Let $l_{i_*(w)}$ be the oriented based loop representing a cyclic geodesic included word $i_*(w)$ which has letters of a single parity in the augmented Dehn presentation of a prime alternating link. Let $S$ be the set of oriented loops given by forgetting the base point of the loops constructed from  $l_{i_*(w)}$ by all possible repeated applications of the method described above.  Then $S$ represents all cyclic geodesic elements of the conjugacy class of $i_*(w)$.
\end{Lem}

Now suppose that $w$ contains letters of both parities,  $n_{\omega}$ letters of parity white and  $n_{\beta}$ letters of parity black.
Since we want to use the algorithm to generate all the geodesic representatives  of the conjugacy class of $w$ (up to cyclic permutation), rather than producing a diagram, will carry out  algorithm~\ref{joconjalg} in the following way:
find all equivalent geodesics of $w$ (using the geodesic completion).  If any of these words have the parity of the first letter opposite to that of the first letter of $w$, then choose one of them and call it $w_{l1}$, say. Similarly, if any of these words have the parity of the last letter opposite to that of the last letter of $w$, then choose one of them and call it $w_{r1}$.  If $w_{l1}$ exists, then cyclically permute it by one letter (so  the first letter becomes the last).  Call this new word $w^{\prime}_{l1}$.
 Generate all of the equivalent geodesics of $w^{\prime}_{l1}$. If any of these words have their first letter of opposite parity of the first letter of $w$, choose one and call it $w_{l2}$.  Repeat this process for as long as is possible or until $\max \{  n_{\beta}, n_{\omega}  \}$ consecutive steps give no new words.
Do an analogous process for $w_{r1}$.

There are three key step in the above process to interpret geometrically: generating the geodesics, recognizing the parity of the first letter and constructing $w^{\prime}_{lj}$ or $w^{\prime}_{rj}$.
The generation of geodesics was discussed earlier. The parity of the first letter is determined by the first intersection point of the loop. Constructing  $w^{\prime}_{lj}$ from the permutation of  $w_{l j-1}$ and the cyclic permutation corresponds to moving the base point along the loop through one intersection point against the orientation.

Putting this together we get the following procedure:
Given $l_w$ minimal, based and oriented. Note the colour of the first and last regions intersected by $l_w$.  Carry out all possible sequences of type 2 deformations with respect to both checker-board hierarchies which fix the base point. If any of these new loops have the first (resp. last) intersection point of the opposite colour as the first (resp. last) as $l_w$ choose one and move the base point forward (resp. backward) one intersection point and repeat this process.  Continue like this for as long as we can or until we stop obtaining new loops.

Since  algorithm~\ref{joconjalg} produces a set of cyclic geodesic words such that every conjugate cyclic geodesic word is a cyclic permutation of one of those produced, we see that by considering the set of all cyclic permutations of these words we get the set of all cyclic geodesic representatives of the conjugacy class.  This obviously does not depend upon which particular cyclic geodesic we started with.  This gives the following lemma.

\begin{Lem}
Let $l_{i_*(w)}$ be the based oriented loop representing a cyclic geodesic included word $i_*(w)$, which has letters of a both parities, in the augmented Dehn presentation of a prime alternating link. Let $S$ be the set of oriented loops given by forgetting the base point of the loops constructed from  $l_{i_*(w)}$ and carry out all possible repeated applications of type 2 deformations.  Then $S$ represents all cyclic geodesic elements of the conjugacy class of $i_*(w)$.
\end{Lem}

Putting all of this together we obtain the following geometric characterization of the conjugacy classes.

\begin{Prop}
Let $L$ be a prime alternating link, $s \in \pi_1(S^3-L) \subset \pi_1(S^3 - (L \cup \mathcal{O})$ and  $l$ be the oriented loop $s$ determines.
Then by using type 1 reduction and type 2 deformations with respect to both hierarchies of $S^3 -L$, $l$ determines a set of loops with the property that every loop is minimal and represents the cyclic geodesics of the conjugacy class of $s$.  Moreover, every loop representing an element in the conjugacy class of $s$ is equivalent to an element of this set by a sequence of type 1 reductions and type 2 deformations.
\end{Prop}

\sgap

Although it should be fairly clear how to construct a polynomial time  normal surface algorithm for the conjugacy problem for prime alternating links, for completeness we outline  one.

\begin{Alg}
Given two oriented loops $l_{w^{\prime}}$ and  $l_{u^{\prime}}$ in the complement of a reduced prime alternating link. Use algorithm~\ref{wordhiealg} to find equivalent minimal loops $l_{w}$ and  $l_{u}$ respectively.
There are two cases.

The first case is when $l_w$ intersects two surfaces  $F_1$ and $F_i$, for some $i \neq 1$.  Use process X to look for and carry out all possible sequences of type 2 deformations (with respect to the checker-board hierarchy of the same colour as the first surface the loop intersects) .  This produces a finite set of minimal loops.
Choose some base point on each of these loops.  If, when traveling round the loop from the base point in the direction of the orientation,  we meet the surfaces in the hierarchy in the same order and direction as for $l_u$ for any choice of base point then $l_{w^{\prime}}$ and  $l_{u^{\prime}}$ are freely homotopic.  Otherwise they are not.

The second case is when $l_w$ intersects surfaces in the hierarchy of one colour only.  Work with the checker-board hierarchy of the opposite colour.
Use process X to carry out the procedure used in lemma~\ref{normsingle} to produce a finite set of minimal loops.
Choose some base point for each of these loops.  If when traveling round the loop from the base point in the direction of the orientation produces we meet surfaces in the hierarchy in the same order and direction as for $l_u$ for any choice of base point then $l_{w^{\prime}}$ and  $l_{u^{\prime}}$ are freely homotopic.  Otherwise they are not.
\end{Alg}

\sgap

\begin{Remark}
By considering annular diagrams (see \cite{LSbook}) one can show, quite unsurprisingly, that two words in the augmented Dehn presentation of an arbitrary link are conjugate if and only if one can be obtained from the other by a finite sequence of free reductions, exchanging pairs and adding a subword which can be freely reduced.
Equivalently, two loops in a link complement are freely homotopic if and only if there is a finite sequence of type 1 reductions and augmentations and type 2 deformations.
Of course our ability to find such a sequence is dependent upon our ability to solve the conjugacy problem.

A similar statement holds for the word problem.

\end{Remark}

\section{Process X and the Petronio Cell Decomposition} \label{app:petronio}

Process X uses the fact that the hierarchy splits the manifold into 3-balls to check for type 1 and 2 moves.  It is easy to see that the realizability of these moves only depends upon the positions of the endpoints of the arcs and how the hierarchy splits the manifold.
More explicitly, the positions of the surfaces of the hierarchy determine a pattern on the boundary of the 3-balls and process X examines the positions and ordering of the intersections of the arc, which we will call the {\em intersection points}, in this pattern.   By considering the Petronio (or pyramid) cell decomposition (see  \cite{{BPbook},{Pe1992}}) of the link complement, we will fully describe process X for non-split alternating links.

\sgap

Let $L$ be a non-split link and $D$ be its canonical projection.  Construct a cell complex by taking as the  0-cells the North and South poles of the crossing balls. For the 1-cells, take the arcs of the link between the poles of the crossing balls and two 1-cells inside each crossing ball where each of these 1-cells has an end point on each of the poles.
For the 2-cells we take two copies of the each of the regions of the link and modify them slightly so that rather than having  arcs lying on the north-south axes of the crossing balls, the arcs lie on the 1-cells of the complex.  We do this in such a way that any intersection of the 2-cells is at one of these 1-cells.
This 2-complex divides $S^3$ into 3-balls and we take these to be the 3-cells of the complex, where the attaching map is the obvious one.

Some of the 3-cells in this complex have non-trivial intersection with the regions of the link.  We call such 3-cells the {\em sandwiched 3-cells}.

What we have obtained is a cell decomposition of $S^3$ such that $L$ is a sub 1-complex and $N(F_1) \cup_{i=1}^n F_i$ is a sub complex.  In addition  the 0-cells and 1-cells determine graphs on the non-sandwiched 3-cells. These are the ones required by process X.

Observe that if, in the above construction, we only add one 1-cell per crossing ball, one 2-cell per region of $D$ and leave out the sandwiched 3-cells we still get the same graph on the boundary of the 3-balls. This gives  a cell decomposition of $S^3$ with $L$ as a sub-complex.
This is called the {\em Petronio cell decomposition}.
So our questions of intersections of arcs with the checker-board hierarchy may be phrased in terms intersections with the 2-cells of this complex. Finally, since this complex has the property that $S^3 - L = (S^3/L)-\{0-\text{cells}\}$ and we are only concerned with points inside the 2-cells, we see that it is sufficient to use the Petronio cell decomposition. We will denote the 3-cells by $B_+$ and $B_-$, and call the graphs on them the {\em boundary graphs}.

\begin{figure}
\centering
\subfigure[]{\raisebox{2.5mm}{\epsfig{file=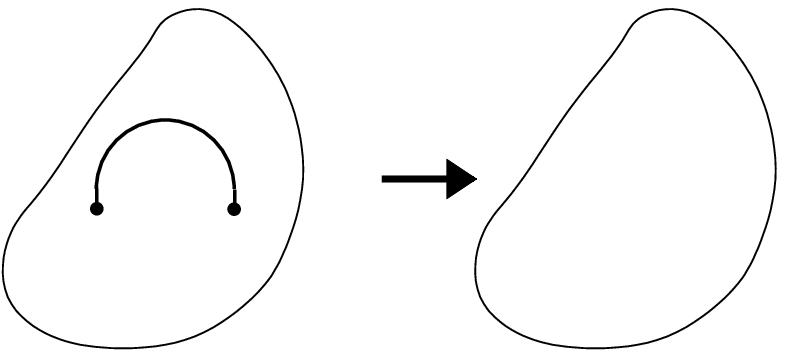, height=1.5cm}}}
\hspace{1cm}
\subfigure[]{\epsfig{file=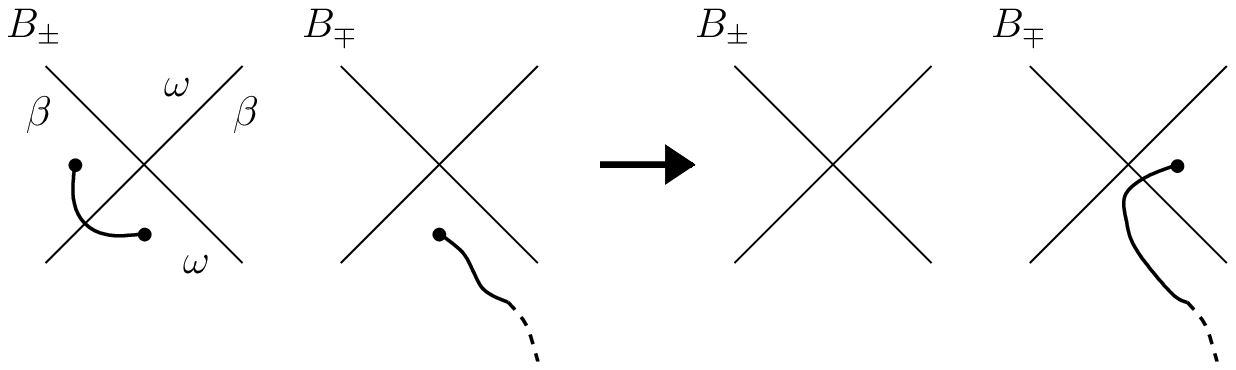, height=2.2cm}}
\caption{The moves on the Petronio cell decomposition.}
\label{bountype}
\end{figure}

\begin{Lem}
Let $L$ be a non-split alternating link.  Then by examining the boundary graphs we see that:

\noindent (a) A type 1 reduction can be recognized by two consecutive intersection points within the same region of the boundary graph. The reduction corresponds to deleting these two points (see figure~\ref{bountype}(a), where the arcs connecting intersection points are to indicate the ordering).

\noindent (b) A type 2 deformation is indicated by two consecutive intersection points in adjacent regions of the boundary path on $B_{\pm}$, which originate from a black region. The deformation replaces these two points with one coming out of the region opposite in the boundary graph on $B_{\mp}$ (see remark~\ref{type2boundarygraph}), as in figure~\ref{bountype}(b).
\end{Lem}

\begin{proof}
The result follows by carrying the processes of  section~\ref{sec:moves} through to the Petronio decomposition.
\end{proof}

\begin{Remark} \label{type2boundarygraph}
When finding the image of a type 2 deformation on a boundary graph we must go back to the link to determine which region we should deform into.  This has no bearing on process X since it only detects the presence of a type 1 or 2 move.
\end{Remark}

\section{Planarity of the Conjugacy Algorithm} \label{app:planar}

In section~\ref{sec:smallconj} it was claimed that the diagrams produced by the conjugacy algorithm for prime alternating links were planar (and therefore the algorithm is polynomial time).  This was proved in \cite{Jo1997}.  Here we present a substantially shorter proof of the planarity of the diagrams produced by algorithm~\ref{joconjalg}.

\medskip

Define the {\em distance} between two regions of a link projection to be the minimum number of edges which a curve between the two regions must cross   ( in terms of the graph of a link with the usual metric it is the length of a geodesic path).

We have seen earlier that  a pair corresponds to a arc of distance 1.
 It is not hard to see that a sister-set corresponds to an arc between two regions which are a distance 2 apart.  Moreover, a region which is distance 1 from each of these regions is the element needed to form the two pairs in the definition of a sister-set.
For completeness we note that a pseudo-pair corresponds to an arc between regions of distance 3.

\begin{The}
The small cancellation diagrams produced by algorithm~\ref{joconjalg} are planar.
\end{The}

\begin{figure}
\centering
\subfigure[]{\raisebox{1.25cm}{\epsfig{file=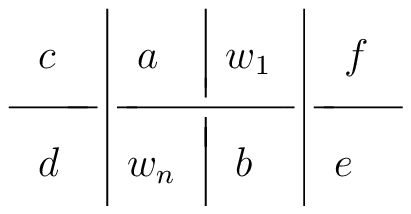, height=1.5cm}}}
\hspace{1cm}
\subfigure[]{\raisebox{5mm}{\epsfig{file=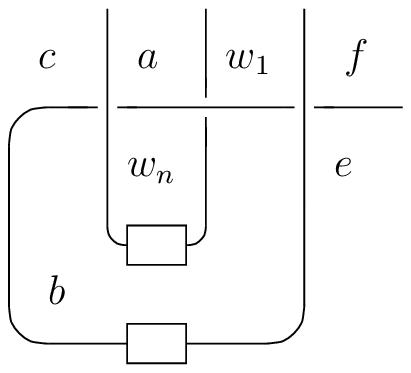, height=3cm}}}
\hspace{1cm}
\subfigure[]{\epsfig{file=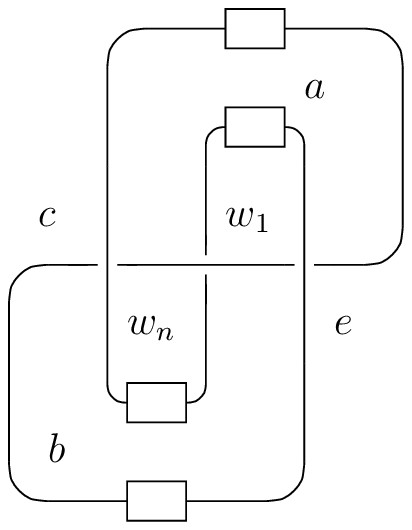, height=4cm}}
\caption{A diagram used in the proof of planarity.}
\label{fig:planarcase1}
\end{figure}

\begin{figure}
\centering
\subfigure[]{\raisebox{1cm}{\epsfig{file=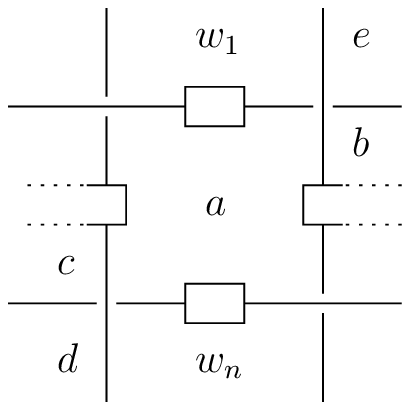, height=4cm}}}
\hspace{1cm}
\subfigure[]{\epsfig{file=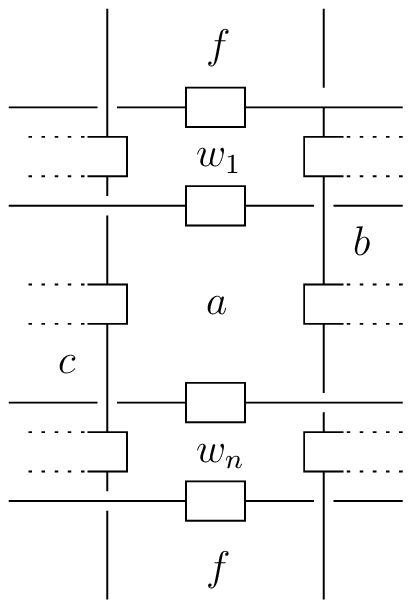, height=6cm}}
\caption{Another diagram used in the proof of planarity.}
\label{fig:planarcase2}
\end{figure}

\begin{proof}
Suppose we are given two cyclic geodesic included words in the \adp of a prime alternating reduced link $L$.  If the words contain letters of both parities the result follows since by theorem~\ref{boundarychains} there are exactly four chains on the boundary (\cite{{Jo1997},{Jo2000}}).  So assume that they both consist of letters of a single  parity. In this case it is sufficient to show that a cyclic geodesic included word  is the inner link path of a conjugacy chain then it is the inner link path of at most two conjugacy chains.

Suppose we are given a  conjugacy chain with inner link path $w=w_1 \cdots w_n$ and chain word $a v a^{-1}$, where $v= v_1 \cdots v_n$. Then
$a^{-1}w_1$ and $w_n^{-1}a$ are pairs and so
$w_i w_{i+1}$ and $v_i v_{i+1}$, where the indices are  modulo $n$, are sister-sets.
Interpreting this into the language of loops  we see that the canonical loop $l_w$  intersects regions of a single colour and consecutive intersection points of the loop occur in regions  a distance of 2 apart.  In particular, the first and last intersection points of $l_w$ occur in regions $w_1$ and $w_n$ a distance of two apart and the region $a$ is adjacent to both of these regions.  So either $a$, $w_1$ and $w_n$ meet at a vertex with $w_1$ and $w_n$ diagonally opposite to each other or they do not meet at a vertex and $a$ is a region which shares edges with $w_1$ and $w_n$.

By  section~\ref{sec:normcon},  we obtain the last and first intersection points of $l_v$ from the last and first of $l_w$ by a type 1 augmentation with respect to region $a$ on an arc which contains the base point, carrying out two type 2 deformations and carrying out a type 1 reduction on the resulting intersection points (which  occurs as the two final intersection points of the new loop). Note that since the link is alternating there are exactly two possible type 2 deformations.
We examine the geometric consequences of these moves.

If the augmented link is the trivial link of two or more components there can be no chains.
There are two cases remaining.

Suppose first that the regions $a$, $w_1$ and $w_n$ meet at a crossing. First consider the case indicated in  figure~\ref{fig:planarcase1}(a).
Since the canonical projection is elementary, $a$, $b$, $w_1$ and $w_n$ are distinct regions.
Then we have $w_n^{-1}w_1 = (w_n^{-1} a)(a^{-1}w_1) = (b^{-1}w_1)(c^{-1}d)$, by hypothesis this implies $b=d$.
This means that there must be a sequence of edges as indicated in figure~\ref{fig:planarcase1}(b), where the boxes indicate the possibility of further knotting.

Now suppose that $w$ is the inner link path of another conjugacy chain. This chain must have its chain word of the form $b u b^{-1}$, for some word $u$, since $a$ and $b$ are the only regions distance one from $w_1$ and $w_n$.
This chain gives $w_n^{-1}w_1 = (w_n^{-1} b)(b^{-1}w_1) = (f^{-1}e)(w_n^{-1}a)$, and hypothesis this implies that $a=f$ and the planar projection is of the form indicated in figure~\ref{fig:planarcase1}(c).
Now since $a \neq b$ it is clear that there can be no more than two conjugacy chains.

The arguments for $w_nw_1^{-1}$, the other positions for $a$ and for the opposite sign of the crossing are similar.

The second case is when the regions $a$, $w_1$ and $w_n$ do not meet at a crossing.  Consider the situation where the projection locally looks like figure~\ref{fig:planarcase2}(a) (note there are necessarily an even number of crossings with the edges of $a$ in the boxes).
Then  $w_n^{-1}w_1 = (w_n^{-1} a)(a^{-1}w_1) = (e^{-1}b)(c^{-1}d)$, which by hypothesis implies that $b=c$.

Now suppose that $w$ is the inner link path of another conjugacy chain. This chain must have its chain word of the form $f u f^{-1}$, for some word $u$. The the region $f$ can only be positioned as in figure~\ref{fig:planarcase2}(b).  This separates the part of the link on the left from the part of the link on the right of the projection as indicated in the figure.  But this is impossible unless $b=c=f$, and by studying the figure it is clear that there are no more suitable regions in which to make a conjugacy chain. Thus there are at most two.

The arguments for $w_nw_1^{-1}$ and the opposite sign of  crossings are similar.
This completes the proof.
\end{proof}

\def\c4t4{$C^{\prime \prime} (4) - T (4)$ }
\def\plane{\mathbb{R}^2}
\def\threespace{\mathbb{R}^3}
\def\doublel{D_i  (L)}

\section{A New Proof that Alternating Links are Non-Trivial}

There are several approaches  in the literature for showing that alternating links are non-trivial - using the determinant in \cite{Ba1930}, the Alexander polynomial in \cite{Cr1959} and \cite{Mu1958}, the Jones polynomial in \cite{Ka}, the Q-polynomial in \cite{Ki1987} and geometrically in \cite{MT1991}.
These proofs give varying degrees of intuition as to why the result holds.
The argument presented here uses Dehn's lemma and the solvability of the word problem to show in a very direct way that spanning disks for the link can not exist giving a very intuitive feel for  the non-triviality of alternating links. We prove:
\begin{The} \label{th:main}
If $L$ is a link admitting a reduced, prime, alternating projection, then $L$ is non-trivial.
\end{The}
Note that since the connected sum of two non-trivial links is non-trivial  restricting ourselves to prime links does not compromise the spirit of this section.

The $i$-th {\em double} of a link $L$ is a parallel copy of the $i$-th component. We define
the $i$-th longitude $\lambda_i$ to be an element of the link group determined by the $i$-th double.
Note that as in Chapter~\ref{chapter:mu}, our longitudes not necessarily null-homologous in the link complement. This is to simplify the argument and causes no real problems.
The following piece of folklore is a consequence of Dehn's lemma and the loop theorem.
\begin{Folklore}
A link is trivial if and only if all of its longitudes are trivial  in the link group.
\end{Folklore}
This reduces theorem~\ref{th:main} to solving the word problem for the  longitudes of the link, which we shall do using small cancellation theory.

\medskip

Recall that the  checker-board colouring of a link projection is an assignment of a  colour black or white to each of the regions of the projection in such a way that  adjacent regions are assigned  a different colour.
\begin{Lem} \label{lem:zigzag}
The $i$-th double of an alternating link $L$ is isotopic to a simple closed curve $J \subset S^3 -L$ such that, in terms of the projection,
any ``intersections'' of $J$ with white regions of the checker-board colouring occur before any ``intersections'' with black regions, with respect to a chosen base point and orientation.
\end{Lem}

\begin{figure}
\centering
\epsfig{file=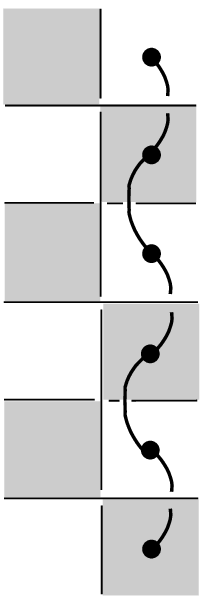, height=4cm}
\hspace{4mm}
\raisebox{2cm}{$\longrightarrow$}
\hspace{4mm}
\epsfig{file=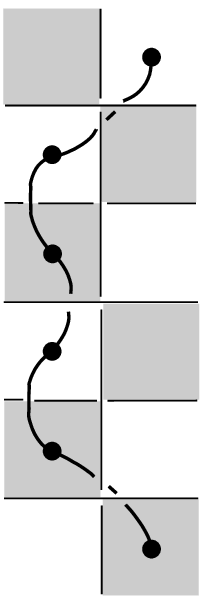, height=4cm}
\hspace{4mm}
\raisebox{2cm}{$\longrightarrow$}
\hspace{4mm}
\epsfig{file=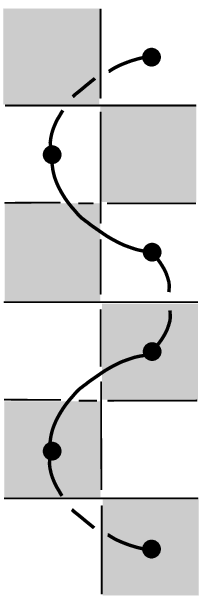, height=4cm}
\caption{The argument of lemma~\ref{lem:zigzag}.}
\label{wobble}
\end{figure}

\begin{proof}
We work in terms of the alternating projection  of the  link $L$ and its $i$-th double $\doublel$. We give the projection of $L$ the checker-board colouring and look at the way that $\doublel$ intersects the regions.

First observe that $\doublel$ is a curve which travels parallel to the component of $L$ and intersects adjacent regions of the projection.
Choose a base point and orientation of $\doublel$ and label the intersection points of $\doublel$ with the regions $v_1, \ldots ,v_n$ of the projection of $L$, where we travel in the direction of the orientation from the base point. We may assume that $v_1$ intersects a white region and therefore $v_n$ intersects a black one.
Since $L$ is alternating we can fix the intersection points $v_1$ and $v_n$ and isotope everything else so that $v_2, \ldots ,v_{n-1}$ all lie in different regions (which are uniquely determined).
Notice that $v_2$ now lies in a white region and $v_{n-1}$ lies in a black region.
Now fix $v_2$ and $v_{n-1}$ and isotope  so that $v_3, \ldots ,v_{n-2}$ lie in different regions (in fact the regions they were originally in). This places $v_3$ in a white region and $v_{n-2}$ in a black region.
Repeating this process a finite number of times gives the required curve $J$.

This argument is embodied in figure~\ref{wobble}.
\end{proof}

We can now prove the main result of this subsection.
\begin{proof}[proof of theorem \ref{th:main}] \footnote{16/11/05: This proof contains a small gap! See my paper {\em A new proof that alternating links are non-trivial} where this gap was plugged by a small extension of the argument.   }
By Lemma~\ref{lem:zigzag} and the geometric interpretation of the generators of the augmented Dehn presentation in section~\ref{sec:loo}, the longitude can be represented by the conjugate of a non-empty word $w$ which changes parity exactly once.
Since the projection is reduced, $w$ is freely reduced.
 A word of this form can not contain a chain word (as these change parity twice) and since the augmented Dehn presentation is a \c4t4 small cancellation group (lemma~\ref{dehncan}), the geodesic characterization theorem tells us that the longitudes are non-trivial.
\end{proof}

\bibliographystyle{plain}

\end{document}